\documentclass[11pt,twoside]{article}%
\usepackage{enumerate}
\usepackage{enumitem}
\usepackage{amsmath,amssymb,amsfonts}
\usepackage{amsfonts}%
\usepackage{graphicx}
\usepackage{color}
\usepackage{pdfsync}
\providecommand{\U}[1]{\protect\rule{.1in}{.1in}}
\newtheorem{thm}{Theorem}[section]
\newtheorem{lem}{Lemma}[section]

\newtheorem{prop}{Proposition}[section]
\newtheorem{corol}{Corollary}[section]
\newtheorem{claim}{Claim}[section]

\newenvironment{pf}[1][\bfseries Proof]{\noindent{#1.} }{\hfill \rule{0.5em}{0.5em}\\}

\numberwithin{equation}{section}
\begin{document}
\title{Existence of multi-peak solutions for a class of quasilinear problems in Orlicz-Sobolev spaces}

\author{Claudianor O. Alves\thanks{C.O. Alves was partially supported by CNPq/Brazil  304036/2013-7  and INCT-MAT, \hspace*{.7cm} e-mail:  coalves@mat.ufcg.edu.br},~~  Ailton R. da Silva\thanks{A.R. da Silva,~ e-mail: ardsmat@gmail.com}\\
Universidade Federal de Campina Grande\\
 Unidade Acad\^emica de Matem\'atica - UAMat\\
 CEP: 58.429-900 - Campina Grande - PB - Brazil}
\date{}
\maketitle

\begin{abstract}
The aim of this work is to establish the existence of multi-peak solutions for the following class of quasilinear problems
\[
- \mbox{div}\big(\epsilon^{2}\phi(\epsilon|\nabla u|)\nabla u\big) + V(x)\phi(\vert u\vert)u  =  f(u)\quad \mbox{in} \quad \mathbb{R}^{N},
\]
where $\epsilon$ is a positive parameter, $ N\geq 2$, $V, f$ are continuous functions satisfying some technical conditions and $\phi$ is a $C^{1}$-function.

\end{abstract}



\section{Introduction}

Several recent studies have focused on the nonlinear Schr\"{o}dinger
equation
$$
i\epsilon\displaystyle \frac{\partial \Psi}{\partial t}=-\epsilon^{2}\Delta
\Psi+(V(z)+E)\Psi-f(\Psi)\,\,\, \mbox{for all}\,\,\, z \in
\mathbb{R}^{N},\eqno{(NLS)}
$$
where $N \geq 2$, $\epsilon > 0$ is a positive parameter and $V,f$ are continuous function verifying some conditions. This class of equation is one of the main objects of the quantum physics, because it appears in problems involving nonlinear optics, plasma physics and condensed matter physics.

The knowledge of the solutions for the elliptic equation
$$
\ \  \left\{
\begin{array}{l}
- \epsilon^{2} \Delta{u} + V(z)u=f(u)
\ \ \mbox{in} \ \ \mathbb{R}^{N},
\\
u \in H^{1}(\mathbb{R}^{N}),
\end{array}
\right.
\eqno{(S)_{\epsilon}}
$$
has a great importance in the study of standing-wave solutions
of $(NLS)$. The existence and concentration of
positive solutions for general semilinear elliptic equations
$(S)_\epsilon$ for the case $N \geq 3$ have been extensively
studied, see for example, Floer and Weinstein \cite{FW}, Oh
\cite{O2}, Rabinowitz \cite{rabinowitz}, Wang \cite{WX}, Cingolani and Lazzo \cite{CL97},
Ambrosetti, Badiale and Cingolani  \cite{ABC}, Floer and Weinstein \cite{FN}, Gui \cite{G}, del Pino and Felmer \cite{DF1}  and their references.

In the above mentioned papers, the existence,  multiplicity and concentration of positive solutions have been obtained in connection with the geometry of the function $V$. In \cite{rabinowitz}, by a mountain pass argument,  Rabinowitz proves the existence of positive solutions of $(S)_{\epsilon}$ for $\epsilon > 0$  
small and
$$
\liminf_{|z| \rightarrow \infty} V(z) > \inf_{z \in
	\mathbb{R}^N}V(z)=V_{0} >0.\eqno{(R)}
$$
Later Wang \cite{WX} showed that these solutions concentrate at global minimum points of $V$ as  $\epsilon$ tends to 0. In \cite{DF1}, del Pino and Felmer have found solutions which concentrate around local minimum of $V$ by introducing a penalization method. More precisely, they assume that there is an open and bounded set $\Lambda \subset \mathbb{R}^N$ such that
$$
0< V_{0}\leq \inf_{z\in\Lambda}V(z)< \min_{z \in
	\partial\Lambda}V(z).
$$

The existence of multi-peak solution has been considered in some papers.
In \cite{G}, Gui has showed the existence of a $\kappa$-peak solution $u_{\epsilon}$ for the problem
$(S)_{\epsilon}$ under the assumptions that $V:\mathbb{R}^{N} \to \mathbb{R}$ is a continuous function verifying
\[
V(z)\geq V_{0}>0\,\,\,\mbox{for all}\,\, z \in \mathbb{R}^{N}\eqno{(V_{0})}
\]
and there exist $\kappa$ disjoint bounded regions $\Omega_{1}, ..., \Omega_{\kappa}$ such that
\[
M_{i} = \min_{x \in \partial \Omega_{i}}V(x) > \alpha_{i} = \inf_{x \in \Omega_{i}}V(x)\,\,\, i= 1,..., \kappa .\eqno{(V_{1})}
\]
A similar result was also obtained by  del Pino and Felmer in \cite{DF2} by using a different approach. In \cite{Alves2011}, Alves has generalized the results found in \cite{G} for a  class of quasilinear problems involving the $p$-Laplacian operator. The reader can find more information about multi-peak solutions for quasilinear problems associated with $(S)_{\epsilon}$ in  Giacomini and Squassina \cite{GS}, Zhang and Xu \cite{Z.Hao} and their references.

After a bibliography review, we did not find any paper related to the existence of multi-peak solution for quasilinear problems involving $N$-functions. Motivated by this fact, we are interesting in finding  multi-peak positive solutions for the following class of quasilinear problems
\begin{align}
\left\{
\begin{array}
[c]{rcl}%
- \mbox{div}\big(\epsilon^{2}\phi(\epsilon|\nabla u|)\nabla u\big) + V(x)\phi(\vert u\vert)u  &=&  f(u)~  \mbox{in}~ \mathbb{R}^{N},\\
u \in W^{1, \Phi}(\mathbb{R}^{N}), &  &
\end{array}
\right. \tag{$ P_{\epsilon} $}\label{PP3}%
\end{align}
where $\epsilon$ is a positive parameter, $ N\geq 2$, $V : \mathbb{R}^{N} \rightarrow \mathbb{R} $ is a continuous function satisfying $(V_{0})$ and $(V_{1})$, and $\phi: [0, +\infty)\rightarrow [0, +\infty)$ is a $C^{1}$-function verifying:
\begin{enumerate}[label={($\phi_\arabic{*}$})]
	\setcounter{enumi}{0}
	\item\label{H1} $\phi(t)$, $(\phi(t)t)^{'}>0$ for all $t>0$.
	
	\item\label{H2} There exist $l,m \in (1, N)$, such that $l\leq m <l^{*}=\frac{Nl}{N-l}$ and
	\[
		l\leq \frac{\phi(t)t^{2}}{\Phi(t)}\leq m \quad \forall t\neq 0,\quad \mbox{where}\,\,\, \Phi(t) = \int_{0}^{|t|}\phi(s)sds.
	\]
	\item\label{H3} The function $\displaystyle\frac{\phi(t)}{t^{m-2}}$ is nonincreasing in $(0, +\infty)$.
	
	\item\label{H4} The function $\phi$ is monotone.
	
	\item\label{H5} There exists a constant $c>0$ such that
	\[
	\vert \phi^{'}(t)t\vert \leq c\phi(t), \quad \forall \ t \in [0, +\infty).
	\]
\end{enumerate}

\noindent Hereafter, we will say that $\Phi \in \mathcal{C}_m$ if
$$
\Phi(t) \geq |t|^{m}, \,\,\,\,\,\, \forall t \in \mathbb{R}. \leqno{(\mathcal{C}_m)}
$$
Moreover, let us denote by $\gamma$ the following real number
$$
\gamma=
\left\{
\begin{array}{l}
m, \,\,\,\, \mbox{if} \,\,\, \Phi \in \mathcal{C}_m, \\
\mbox{}\\
l, \,\,\,\, \mbox{if} \,\,\, \Phi \notin \mathcal{C}_m.
\end{array}
\right.
$$
Related to the function $f$, we assume that it is a $C^{1}$-function satisfying
\begin{enumerate}[label={($f_\arabic{*}$})]
	\item\label{f1}  There are functions $r,b:[0,+\infty) \to [0,+\infty)$ such that
	\[
	\limsup_{|t|\rightarrow 0} \frac{f^{'}(t)}{(r(|t|)|t|)^{'}}=0 \quad \mbox{and} \quad \limsup_{|t|\rightarrow +\infty} \frac{|f^{'}(t)|}{(b(|t|)|t|)^{'}}<+\infty.
	\]
	\item\label{f2} There exists $\theta > m$ such that
	\[
	0< \theta F(t) \leq f(t)t \quad \forall t>0,\quad \mbox{where}\quad   F(t)= \int_{0}^{t}f(s)ds.
	\]
	\item\label{f3}  The function $\displaystyle\frac{f(t)}{t^{m -1}}$ is increasing for $t>0$.
\end{enumerate}
Here, the functions $r $ and $b$ are $C^{1}$-function satisfying the following conditions:
\begin{enumerate}[label={($r_\arabic{*}$})]
	\setcounter{enumi}{0}
	\item\label{R1} $r$ is increasing.
	\item\label{R2} There exists a constant $\overline{c}>0$ such that
	\[
	\vert r^{'}(t)t \vert \leq \overline{c} r(t), \quad \forall \ t\geq 0.
	\]
	\item\label{R3} There exist positive constants $r_{1}$ and $r_{2}$ such that
	\[
	r_{1}\leq \frac{r(t)t^{2}}{R(t)}\leq r_{2}, \quad \forall t\neq 0,\quad \mbox{where}\quad R(t) = \int_{0}^{|t|} r(s)sds.
	\]
	\item\label{R4} The function $R$ satisfies
	\[
	\limsup_{t\rightarrow 0}\frac{R(t)}{\Phi(t)}< +\infty \quad \mbox{and} \quad \limsup_{|t|\rightarrow +\infty}\frac{R(t)}{\Phi_{*}(t)}=0.
	\]
	
\end{enumerate}

\begin{enumerate}[label={($b_\arabic{*}$})]
	\setcounter{enumi}{0}
	\item\label{B1} $b$ is increasing.
	\item\label{B2} There exists a constant $\widetilde{c}>0$ such that
	\[
	\vert b^{'}(t)t\vert \leq \widetilde{c}b(t), \quad \forall t\geq 0.
	\]
	\item\label{B3} There exist positive constants $b_{1}, b_{2} \in (1, \gamma^*)$ such that
	\[
	b_{1}\leq \frac{b(t)t^{2}}{B(t)}\leq b_{2} \;\; \forall t\neq 0,\;\; \mbox{where}\;\; B(t) = \int_{0}^{|t|} b(s)sds\;\; \mbox{and}\;\; \gamma^*=\frac{N\gamma}{N-\gamma}.
	\]
	\item\label{B4} The function $B$ satisfies
	\[
	\limsup_{t\rightarrow 0}\frac{B(t)}{\Phi(t)}< +\infty \quad \mbox{and} \quad \limsup_{|t|\rightarrow +\infty}\frac{B(t)}{\Phi_{*}(t)}=0,
	\]
\end{enumerate}
where $\Phi_{*}$ is the Sobolev conjugate function, which is defined as being the inverse function of
\[
G_{\Phi}(t) = \int_{0}^{t}\frac{\Phi^{-1}(s)}{s^{1+\frac{1}{N}}}ds.
\]

Using the change variable $v(x)=u({x}/{\epsilon})$, it is easy to see that the problem \eqref{PP3} is equivalent to the following problem
\begin{align}
\left\{
\begin{array}
[c]{rcl}%
- \Delta_{\Phi}v + V(\epsilon x)\phi(|v|)v & = & f(v)~  \mbox{in}~ \mathbb{R}^{N},\\
v \in W^{1, \Phi}(\mathbb{R}^{N}), &  &
\end{array}
\right. \tag{$ \widetilde{P}_{\epsilon} $}\label{pp3}%
\end{align}
where the operator $ \Delta_{\Phi} u  = \mathrm{div}(\phi(\vert\nabla u\vert)\nabla u) $, named $\Phi$-Laplacian operator, is a natural extension of the $p$-Laplacian operator, with $ p $ being a positive constant. This operator appears in a lot of physical applications, such as

\vspace{.3cm}
\noindent {\it Non-Newtonian Fluid:  } \, $\Phi(t) = \frac{1}{p}|t|^{p}$ for $p>1$,\\
\noindent { \it Plasma Physics: } \, $\Phi(t) = \frac{1}{p}|t|^{p} + \frac{1}{q}|t|^{q}$ where $1<p<q<N$ with $q \in (p, p^{*})$,
\noindent {\it Nonlinear Elasticity:} \, $ \Phi(t) = (1+t^{2})^{\alpha}-1, \alpha \in (1, \frac{N}{N-2})$, \\
\noindent { \it Plasticity:} \, $ \Phi(t) = t^{p}\ln(1+t), 1< \frac{-1+\sqrt{1+4N}}{2}<p<N-1, N\geq 3$, \\
\noindent { \it Generalized Newtonian Fluid:} \, $\Phi(t) = \int_{0}^{t}s^{1-\alpha}(\sinh^{-1}s)^{\beta}ds, 0\leq \alpha \leq 1$ and $\beta > 0$.

\vspace{.3cm}
The reader can find more details about the physical applications in \cite{Db}, \cite{FN}, \cite{FN2} and their references. The existence of solution for \eqref{pp3} when $\epsilon=1$ in bounded and unbounded domains of $\mathbb{R}^{N}$ has been established in some paper, see for example \cite{AGJ}, \cite{AdP}, \cite{BBR1}, \cite{BBR2}, \cite{FN1}, \cite{LMM}, \cite{MRR}, \cite{MR1}, \cite{MR2}, \cite{J}, \cite{JS} and references therein. However, associated with the existence, multiplicity and concentration of solution for a $\Phi$-Laplacian equation, the authors know only the papers \cite{AS} and \cite{AS2}.

\vspace{0.5 cm}

Now, we are ready to state our main result.

\begin{thm}\label{T1}
Suppose that  \ref{H1}-\ref{H3}, \ref{B1}-\ref{B3}, \ref{f1}-\ref{f3}, $(V_{0})$ and $(V_{1})$ hold. Then, for each $\Gamma \subset \big\{ 1, ..., \kappa\big\}$, there exist $\epsilon^{*}>0$ such that, for that, for $ \epsilon \in (0, \epsilon^{*}]$, \eqref{PP3} has a family $\{u_{\epsilon}\}$ of positive solutions verifying the following property for $\epsilon$ small enough:\\
There exists $\delta > 0$ such that
\[
\sup_{x \in \mathbb{R}^{N}}u_{\epsilon}(x)\geq \delta.
\]
There exists $P_{\epsilon, i} \in \Omega_{i}$ for all $i \in \Gamma$ such that, for each $\eta > 0$, there exists $\rho>0$ verifying
\[
\sup_{x \in B_{\epsilon \rho}(P_{\epsilon, i})}u_{\epsilon}(x)\geq \delta\,\,\, \mbox{for all}\,\, i \in \Gamma
\]
and
\[
\sup_{x \in \mathbb{R}^{N}\backslash \cup_{i \in \Gamma}B_{\epsilon \rho}(P_{\epsilon, i})}u_{\epsilon}(x)< \eta.
\]
\end{thm}

In the above theorem, if $\Gamma$ has $\iota$ elements, we say that $u_{\epsilon}$ is a $\iota$-peak solution.  From now on, we will work with \eqref{pp3} to get multi-peak solutions of \eqref{PP3}.

The proof of our main theorem will make by using variational methods and adpating some arguments found in   \cite{Alves2011}  and \cite{G}. However, we would like point out that some estimates in our paper are totally different from those used in \cite{Alves2011, G}, because some properties and estimates that occur for the $p$-Laplacian do not hold for a $\Phi$-Laplacian equation. Here, we overcome these difficulties by showing a new version of Lions' Lemma for Orlicz-Sobolev Spaces and also a new property involving the Orlicz-Sobolev, these two results can be seen in Section 5. Moreover, in \cite{Alves2011} was used the interaction Moser techniques, which does not work well in our case. Hence, it was necessary to change the arguments and we have used some ideas found in  \cite{Fusco} and \cite{LU}.

Before concluding this section, we would like to say that the reader can find a brief review about Orlicz-Sobolev spaces in \cite{AS}, \cite{AS2} and \cite{FN1}. However, for a more detailed study, we cite the books \cite{Adams1}, \cite{MU} and \cite{Rao}.

\vspace{0.5 cm}

\noindent \textbf{Notation:} In this paper, we use the following
notations:
\begin{itemize}
	\item  If $A$ is a N-function, we denote by $\tilde{A}$ and $A_*$ its complementary and conjugate functions respectively.
	
	\item  If $A$ is a N-function, we denote by $L^{A}(\mathbb{R}^N)$ and $W^{1,A}(\mathbb{R}^N)$ the Orlicz and Orlicz-Sobolev spaces respectively. Moreover, we denote by $\|\,\,\,\,\,\|_{A}$ and $\|\,\,\,\,\,\|_{1,A}$ their usual norms given by 
$$
\|  u \|_{A}:= \inf\Big\{ \lambda > 0;  \int_{\mathbb{R}^{N}}A\Big(\frac{\vert  u\vert}{\lambda}\Big)dx \leq 1\Big\}
$$
and
$$
\|u\|_{1,A}=\|\nabla u\|_A+\|u\|_A.
$$
\item We say that a N-function $A$ verifies the $\Delta_{2}$-condition, denote by  $A \in \Delta_{2}$, if there is a $K > 0$ such that
\[
A(2t) \leq K A(t),\quad \forall t\geq 0.
\]
If $A, \tilde{A}$ are  N-functions verifying the $\Delta_{2}$-condition, then $L^{A}(\mathbb{R}^N)$ and $W^{1,A}(\mathbb{R}^N)$  are reflexive and separable. As $(\phi_1)$-$(\phi_2)$ imply that $\Phi$ and $\tilde{\Phi}$ satisfy the $\Delta_{2}$-condition, we have $L^{\Phi}(\mathbb{R}^N)$ and $W^{1,\Phi}(\mathbb{R}^N)$  are reflexive and separable spaces.

\item   $C$ denotes (possible different) any positive constant, whose value is not relevant.

\end{itemize}

\section{Penalization Method}
In present section our main goal is to prove the existence of solution for an auxiliary problem by adapting some ideas explored in \cite{Alves2011} and \cite{G}.

Since we intend to find positive solutions, we will assume that
\begin{eqnarray}\label{posf2}
f(t) = 0, \quad \forall t \in (-\infty, 0].
\end{eqnarray}

In what follows, let us denote by $I_{\epsilon} : X_{\epsilon} \rightarrow \mathbb{R}$ the energy functional  given by
\[
I_{\epsilon}(u)= \int_{\mathbb{R}^{N}}\Phi(\vert\nabla u\vert)dx + \int_{\mathbb{R}^{N}}V(\epsilon x)\Phi(\vert u\vert)dx - \int_{\mathbb{R}^{N}}F(u)dx,
\]
where $X_{\epsilon}$ denotes the subspace of $W^{1, \Phi}(\mathbb{R}^{N})$ given by
\[
X_{\epsilon}= \Big\{ u \in W^{1, \Phi}(\mathbb{R}^{N})\, :\, \int_{\mathbb{R}^{N}}V(\epsilon x)\Phi(\vert u \vert)dx < +\infty\Big\}
\]
endowed with the norm
\[
\Vert u\Vert_{\epsilon} = \Vert \nabla u\Vert_{\Phi} + \Vert u\Vert_{\Phi, V_{\epsilon}}
\]
where
$$
\Vert \nabla u \Vert_{\Phi}:= \inf\Big\{ \lambda > 0;  \int_{\mathbb{R}^{N}}\Phi\Big(\frac{\vert \nabla u\vert}{\lambda}\Big)dx \leq 1\Big\}
$$
and
\[
\Vert u\Vert_{\Phi, V_{\epsilon}}:= \inf\Big\{ \lambda > 0;  \int_{\mathbb{R}^{N}}V(\epsilon x)\Phi\Big(\frac{\vert u\vert}{\lambda}\Big)dx \leq 1\Big\}.
\]
As $\Phi$ and $\widetilde{\Phi}$ verify $\Delta_2$-condition, the space $X_\epsilon$ is reflexive and separable. Moreover, from $(V_{0})$, it follows that the embeddings
\begin{eqnarray}\label{IC}
X_{\epsilon}\hookrightarrow L^{\Phi}( \mathbb{R}^{N}) \quad \mbox{and}\quad X_{\epsilon}\hookrightarrow L^{B}( \mathbb{R}^{N})
\end{eqnarray}
are continuous. From the above embeddings, a direct computation yields   $I_{\epsilon} \in C^{1}(X_{\epsilon}, \mathbb{R})$ with
\[
I^{'}_{\epsilon}(u)v = \int_{\mathbb{R}^{N}}\phi(\vert\nabla u\vert)\nabla u \nabla v \, dx + \int_{\mathbb{R}^{N}}V(\epsilon x)\phi(\vert u\vert)uv \, dx - \int_{\mathbb{R}^{N}}f(u)v \, dx
\]
for all $u,v \in X_{\epsilon}$. Thereby, $u \in X_{\epsilon}$ is a weak solution of \eqref{pp3} if, and only if, $u$ is a critical point of $I_{\epsilon}$. Furthermore, by \eqref{posf2}, the critical points of $I_{\epsilon}$ are
nonnegative.

Let $\theta$ be the number given in \ref{f3} and $a, \xi > 0$ satisfying
$$
\xi >\displaystyle\frac{(\theta - l)}{(\theta - m)}\frac{m}{l} \quad \mbox{and} \quad \displaystyle\frac{f(a)}{\phi(a)a} = \frac{V_{0}}{\xi}.
$$
Using the above numbers, let us define the function
\begin{align}
\widetilde{f}(s) = \left\{
\begin{array}
[c]{rcl}%
f(s) && if ~  s\leq a,\\
\displaystyle\frac{V_{0}}{\xi}\phi(s)s  && if ~ s>a.
\end{array}\nonumber
\right.%
\end{align}
Finally, fixed $\Gamma \subset \big\{1, ..., \kappa\big\}$, we consider the function 
\[
g(x, s) = \chi_{\Omega}(x)f(s) + (1-\chi_{\Omega}(x))\widetilde{f}(s),
\]
where $\chi_{\Omega}$ is the characteristic function related to the set 
\[
\Omega = \bigcup_{i \in \Gamma}\Omega_{i}.
\]

 From definition of $g$, it follows that  $g$ is a Carath\'eodory function verifying
\begin{eqnarray}\label{posg3}
g(x, s) = 0, \quad \forall (x,s) \in \mathbb{R}^{N} \times (-\infty, 0]
\end{eqnarray}
and
\begin{eqnarray}\label{gf3}
g(x, s) \leq f(s), \quad \forall (x,s) \in \mathbb{R}^{N} \times \mathbb{R}.
\end{eqnarray}
Moreover,  the following conditions also hold:
\begin{enumerate}[label={($g_\arabic{*}$})]
	\setcounter{enumi}{0}
	
	\item\label{g13}$0\leq\theta G(x, s) = \theta\displaystyle\int_{0}^{s}g(x, t)dt \leq g(x, s)s$, $\forall (x,s) \in \Omega \times (0,+\infty)$.
	\item\label{g23} $0<lG(x, s)\leq g(x, s)s\leq \displaystyle\frac{V_0}{\xi}\phi(s)s^{2}$, $\forall (x,s) \in \Omega^{c} \times (0,+\infty)$.
	
\end{enumerate}

Using the function $g$, we set the auxiliary problem
\begin{align}
\left\{
\begin{array}
[c]{rcl}%
- \Delta_{\Phi}u + V(\epsilon x)\phi(|u|)u & = & g(\epsilon x, u)~  \mbox{in}~ \mathbb{R}^{N},\\
u \in W^{1, \Phi}(\mathbb{R}^{N}). &  &
\end{array}
\right. \tag{$ {A}_{\epsilon} $}\label{PA3}%
\end{align}
Associated with  \eqref{PA3}, we have the functional $J_{\epsilon}: X_{\epsilon}\rightarrow \mathbb{R}^{N}$ given by
\[
J_{\epsilon}(u)= \int_{\mathbb{R}^{N}}\Phi(\vert\nabla u\vert)dx + \int_{\mathbb{R}^{N}}V(\epsilon x)\Phi(\vert u\vert)dx - \int_{\mathbb{R}^{N}}G(\epsilon x, u)dx,
\]
which belongs to $C^{1}(X_{\epsilon}, \mathbb{R})$ with
\[
J^{'}_{\epsilon}(u)v = \int_{\mathbb{R}^{N}}\phi(\vert\nabla u\vert)\nabla u \nabla v \, dx + \int_{\mathbb{R}^{N}}V(\epsilon x)\phi(\vert u\vert)uv \, dx - \int_{\mathbb{R}^{N}}g(\epsilon x, u)v \, dx,
\]
for all $u,v \in X_{\epsilon}$. Therefore, critical points of $J_{\epsilon}$ are nonnegative weak solutions of \eqref{PA3}.

Here, we would like to point out that if $u_{\epsilon}$ is a positive solution of \eqref{PA3} with $u_{\epsilon}(x)\leq a$ for every $x \in \mathbb{R}^{N}\backslash \Omega_{\epsilon}$ with $\Omega_{\epsilon} = \Omega / \epsilon$, then $u_{\epsilon}$ is also a positive solution of \eqref{pp3}.

\subsection{The behavior of the $(PS)^{*}_{c}$ sequences}

In what follows, we say that $(u_{n})$ is a $(PS)^{*}_{c}$ sequence when
\[
(u_{n}) \subset X_{\epsilon_{n}},\,\,\, J_{\epsilon_{n}}(u_{n})\rightarrow c,\,\,\,\, \Vert J^{'}_{\epsilon_{n}}(u_{n})\Vert^{*}_{\epsilon_{n}}\rightarrow 0 \,\,\,\, \mbox{and}\,\, \epsilon_{n}\rightarrow 0. \,\,\,
\]
The main result of this section is as follows:
\begin{prop}\label{PropI3}
Let $(u_{n})$ be a $(PS)^{*}_{c}$ sequence. Then, there exist a subsequence of $(u_{n})$, still denoted by itself, a nonnegative integer $p$, sequences of points $(y_{n, j})\subset \mathbb{R}^{N}$ with $j = 1, ..., p$ such that
\[
\epsilon_{n}y_{n, j}\rightarrow x_{j} \in \overline{\Omega}\quad \mbox{and}\quad \vert y_{n, j} - y_{n, i}\vert \rightarrow +\infty\quad \mbox{as} \,\,\,n\rightarrow +\infty
\]
and
\[
\Big\Vert u_{n}(\cdot)-\sum_{j=1}^{p}u_{0, j}(\cdot - y_{n, j})\varphi_{\epsilon_{n}}(\cdot - y_{n, j})\Big\Vert_{\epsilon_{n}} \rightarrow 0\,\,\,\mbox{as}\,\,n \rightarrow +\infty
\]
where $\varphi_{\epsilon}(x) = \varphi\big(x/ (-\ln\epsilon)\big)$ for $0 < \epsilon < 1$, and $\varphi$ is a cut-off function which $\varphi(z) = 1$ for $|z|\leq 1$, $\varphi(z) = 0$ for $|z|\geq 2$ and $\vert \nabla \varphi\vert \leq 2$. The function $u_{0, j} \neq 0$ is a nonnegative solution for
\[
- \Delta_{\Phi}u + V_{j}\phi(\vert u\vert)u  =  g_{0, j}(x, u)\quad \mbox{in} \quad \mathbb{R}^{N},\eqno{(P^{j})}
\]
where $V_{j} = V(x_{j})\geq V_{0}>0$ and $g_{0, j}(x, u) = \displaystyle\lim_{n \rightarrow \infty}g(\epsilon_{n}x + \epsilon_{n}y_{n, j}, u)$. Moreover, we have $c\geq 0$ and
\[
c = \sum_{j = 1}^{p}J_{0, j}(u_{0, j}),
\]
where $J_{0, j} : W^{1, \Phi}(\mathbb{R}^{N}) \rightarrow \mathbb{R}$ denotes the functional given by
\[
J_{0, j}(u) = \int_{\mathbb{R}^{N}}\Phi(\vert\nabla u\vert)dx + V_{j}\int_{\mathbb{R}^{N}}\Phi(\vert u\vert)dx - \int_{\mathbb{R}^{N}}G_{0, j}(x, u)dx
\]
with $G_{0, j}(x, t) = \int_{0}^{t}g_{0, j}(x, s)ds$.
\end{prop}
\begin{pf}
Let $(u_{n})$ be a $(PS)^{*}_{c}$ sequence. Arguing as \cite{AS}, there exists $M>0$ independent of $n$ such that
\begin{eqnarray*}
\Vert u_{n} \Vert_{\epsilon_{n}}\leq M \quad \forall n \in \mathbb{N},
\end{eqnarray*}
showing that $(u_{n})$ is a bounded sequence in $W^{1, \Phi}(\mathbb{R}^{N})$. Since
$$
c+o_n(1)= J_{\epsilon_{n}}(u_{n}) - \frac{1}{\theta}J_{\epsilon_{n}}^{'}(u_{n})u_{n}
$$
and $(\phi_2)$ combined with $(g_1)$-$(g_2)$ give
\begin{eqnarray}\label{jpl32}
J_{\epsilon_{n}}(u_{n}) - \frac{1}{\theta}J_{\epsilon_{n}}^{'}(u_{n})u_{n}\geq C\Bigg(\int_{\mathbb{R}^{N}}\Phi(|\nabla u_{n}|)dx + \int_{\mathbb{R}^{N}}V(\epsilon_{n} x)\Phi(|u_{n}|)dx\Bigg),\quad
\end{eqnarray}
with $C = \Bigg[\Big(1-\frac{m}{\theta}\Big)-\Big(1-\frac{l}{\theta}\Big)\frac{m}{kl}\Bigg]>0$, we deduce that $c\geq 0$. Thereby, if $c = 0$ , \eqref{jpl32} ensures that
\[
\int_{\mathbb{R}^{N}}\Phi(|\nabla u_{n}|)dx \rightarrow 0\quad \mbox{and}\quad \int_{\mathbb{R}^{N}}V(\epsilon_{n} x)\Phi(|u_{n}|)dx\rightarrow 0,
\]
leading to $\Vert u_{n} \Vert_{\epsilon_{n}}\rightarrow 0$. In the sequel, we will consider only the case $c>0$.

We claim that there exist positive constants $\rho, a$, a subequence of $(u_{n})$, still denoted by itself, and a sequence $(y_{n, 1})\subset \mathbb{R}^{N}$ such that
\begin{eqnarray}\label{solposi}
\int_{B_{\rho}(y_{n, 1})}\Phi(\vert u_{n}(x)\vert) dx \geq a>0, \quad \forall n \in \mathbb{N}.
\end{eqnarray}
Otherwise, since $(u_{n})$ is bounded in $W^{1, \Phi}(\mathbb{R}^{N})$, a Lions-type result for Orlicz-Sobolev spaces found in \cite[Theorem 1.3]{AGJ} gives $u_{n}\rightarrow 0$ in $L^{B}(\mathbb{R}^{N})$, that is,
\[
\int_{\mathbb{R}^{N}}B(|u_{n}|)dx \rightarrow 0.
\]
Now, the definition of $J_{\epsilon_{n}}^{'}(u_{n})u_{n}$ together with \ref{f1}, \ref{H2} and \ref{B2} yields 
\[
c_{1}\int_{\mathbb{R}^{N}}B(|u_{n}|)dx + J_{\epsilon_{n}}^{'}(u_{n})u_{n} \geq l\Big(1-\frac{1}{k}\Big)\Big(\int_{\mathbb{R}^{N}}\Phi(|\nabla u_{n}|)dx + \int_{\mathbb{R}^{N}}V(\epsilon_{n} x)\Phi(|u_{n}|)dx\Big),
\]
showing that $\Vert u_{n} \Vert_{\epsilon_{n}}\rightarrow 0$. Then $c=0$, which is a contradiction. Therefore \eqref{solposi} holds.

Now, setting $w_{n, 1}(x)=u_{n}(x + y_{n,1})$, we see that $(w_{n, 1})$ is a bounded sequence in $W^{1, \Phi}(\mathbb{R}^{N})$. Thus, there exist $u_{0, 1} \in W^{1, \Phi}(\mathbb{R}^{N})$ and a subsequence of $(w_{n,1})$, still denoted by itself, such that
\[
w_{n, 1}\rightharpoonup u_{0, 1} \quad \mbox{in}\,\,W^{1, \Phi}(\mathbb{R}^{N}).
\]
The above limit and \eqref{solposi} combine to give $u_{0, 1} \neq 0$.

Hereafter, we will show that $u_{0, 1}$ is the solution of $(P^{1})$. For this purpose, it is crucial to show the following claim:

\begin{claim}\label{Af1}
The sequence $(\epsilon_{n}y_{n, 1})$ is bounded. Moreover, there exists $x_{1} \in \overline{\Omega}$ such that, up to a
subsequence, $\epsilon_{n}y_{n, 1} \rightarrow x_{1}$.
\end{claim}

In fact, suppose by contradiction that $(\epsilon_{n}y_{n, 1})$ is an unbounded sequence. Then, without loss of generality, we can suppose that $|\epsilon_{n}y_{n, 1}|\rightarrow + \infty$. By using the limit  $\displaystyle\lim_{\epsilon \rightarrow 0}\epsilon \ln \epsilon = 0$, it is easy to check that for $n$ large enough
$$
\epsilon_{n}y_{n, 1} + \epsilon_{n}x \in \mathbb{R}^{N} \setminus \Omega \quad \mbox{for} \quad |x|< 2|\ln \epsilon_{n}|.
$$
Once $(u_{n})$ is $(PS)^{*}_{c}$, setting $v_{n}(x) = u_{n}(x)\varphi_{\epsilon_{n}}(x - y_{n, 1})$, we have $(\Vert v_{n}\Vert_{\epsilon_{n}})$ is bounded in $\mathbb{R}$ and $J^{'}_{\epsilon_{n}}(u_{n})v_{n} = o_{n}(1)$. On the other hand, a direct computation gives
$$
J^{'}_{\epsilon_{n}}(u_{n})v_{n} \geq \Big(l - \frac{m}{k}\Big)\int_{\mathbb{R}^{N}}\Big(\Phi(|\nabla u_{0, 1}|) + V_{0}\Phi(| u_{0, 1}|) \Big)dx + o_{n}(1),
$$
and so,
\[
\int_{\mathbb{R}^{N}}\Big(\Phi(|\nabla u_{0, 1}|) + V_{0}\Phi(| u_{0, 1}|) \Big)dx = 0,
\]
implying that $u_{0, 1} = 0$, which is absurd. Thereby, $(\epsilon_{n}y_{n ,1})$ is a bounded sequence. From this, there exists $x_{1} \in \mathbb{R}^{N}$ such that for some subsequence, $\epsilon_{n}y_{n, 1} \rightarrow x_{1}$. The same type of argument works to prove that $x_{1} \in \overline{\Omega}$, which proves the Claim \ref{Af1}.

The same arguments explored in  \cite[Lemma 4.3]{AGJ} work to show that there exists a subsequence of $(w_{n, 1})$, still denote by itself, such that
\begin{eqnarray}\label{DER3}
w_{n, 1}(x)\to u_{0, 1}(x)\,\,\, \mbox{and} \,\,\, \nabla w_{n, 1}(x)\to \nabla u_{0, 1}(x)\,\,\, \mbox{a. e. in}\,\,\mathbb{R}^{N}.
\end{eqnarray}
The Claim \ref{Af1} combined with the limit above permit to conclude that $u_{0, 1}$ is solution of $(P^{1})$.
%

Next, we consider $u_{n}^{1}(x) = u_{n}(x) - (u_{0, 1}\varphi_{\epsilon_{n}})(x - y_{n, 1})$. We will show that $(u_{n}^{1})$ is a $(PS)_{c-J_{0, 1}(u_{0, 1})}^{*}$ sequence, that is,
\begin{eqnarray*}
	J_{\epsilon_{n}}(u_{n}^{1})\rightarrow c - J_{0, 1}(u_{0, 1}) \quad \mbox{and}\quad \Vert J^{'}_{\epsilon_{n}}(u_{n}^{1})\Vert_{\epsilon_{n}}^{*}\rightarrow 0.
\end{eqnarray*}
Firstly, we prove that $J_{\epsilon_{n}}(u_{n}^{1})\rightarrow c - J_{0, 1}(u_{0, 1})$. For this end, from a result due to Brezis and Lieb \cite{BL}, we derive that
\begin{eqnarray}\label{DER4}
J_{\epsilon_{n}}(u_{n}^{1}) - J_{\epsilon_{n}}(u_{n}) + J_{\epsilon_{n}}\big((u_{0, 1}\varphi_{\epsilon_{n}})(x - y_{n, 1})\big) = I_{n} + o_{n}(1),
\end{eqnarray}
where 
$$
I_{n} = \displaystyle\int_{\mathbb{R}^{N}}\Big[ G(\epsilon_{n}x, u_{n}) -  G(\epsilon_{n}x, u^{1}_{n}) -  G(\epsilon_{n}x, (u_{0, 1}\varphi_{\epsilon_{n}})(x - y_{n, 1}))\Big]dx.
$$
Arguing as in \cite[Proposition 2.4]{Alves2011}, given $\eta>0$, there exists $\rho>0$ and $n_{0} \in \mathbb{N}$ such that for $n \geq n_{0}$
$$
	I_{n} \leq  \eta + \int_{|x|\geq \rho}\Big|G(\epsilon_{n}x + \epsilon_{n}y_{n, 1}, w_{n, 1}) - G(\epsilon_{n}x + \epsilon_{n}y_{n, 1}, w_{n, 1}- u_{0, 1}\varphi_{\epsilon_{n}})\Big|dx .
$$
Moreover, increasing $\rho$ if necessary, the conditions \ref{H1}-\ref{H5}, \ref{R1}-\ref{R3}, \ref{B1}-\ref{B3}, \ref{f1}, \ref{f3} and \eqref{gf3} combine to give 
\begin{eqnarray*}
\int_{|x|\geq \rho}\Big|G(\epsilon_{n}x + \epsilon_{n}y_{n, 1}, w_{n, 1}) - G(\epsilon_{n}x + \epsilon_{n}y_{n, 1}, w_{n, 1}- u_{0, 1}\varphi_{\epsilon_{n}})\Big|dx &\leq& \eta + o_{n}(1).
\end{eqnarray*}
On the other hand, a direct calculus given us
\begin{eqnarray}\label{DER5}
J_{\epsilon_{n}}\big((u_{0, 1}\varphi_{\epsilon_{n}})(\cdot - y_{n, 1})\big) \to J_{0, 1}(u_{0, 1}).
\end{eqnarray}
The last inequalities together with \eqref{DER4} and \eqref{DER5} leads to $J_{\epsilon_{n}}(u_{n}^{1})\rightarrow c - J_{0, 1}(u_{0, 1})$. A similar argument can be used to show that $\Vert J^{'}_{\epsilon_{n}}(u_{n}^{1})\Vert_{\epsilon_{n}}^{*}\rightarrow 0$. 

As $(u_{n}^{1})$ is a $(PS)_{c-J_{0, 1}(u_{0, 1})}^{*}$, we can repeat the previous arguments to find a  sequence $(y_{n, 2})\subset \mathbb{R}^{N}$ verifying
\begin{eqnarray}\label{solposi2}
\int_{B_{\rho}(y_{n, 2})}\Phi(\vert u_{n}^{1}(x)\vert) dx \geq a_{1}>0.
\end{eqnarray}
We observe that the sequence $(y_{n, 2})$ can be chosen so that
\begin{eqnarray}\label{DESs3}
\vert y_{n, 2} - y_{n, 1}\vert \to + \infty.
\end{eqnarray}
Indeed, to see why, we assume that $(\vert y_{n, 2} - y_{n, 1}\vert)$ is bounded in $\mathbb{R}$. Thus, by \eqref{solposi2}, there exists $\rho_{1}>0$ such that
$$
\int_{B_{\rho}(y_{n, 2})}\Phi(\vert u_{n}^{1}(x)\vert) dx \leq \int_{B_{\rho_1}(0)}\Phi(|w_{n}(x) - (u_{0, 1}\varphi_{\epsilon_{n}})(x)|)dx
$$
and so,
$$
\int_{B_{\rho_1}(0)}\Phi(|w_{n}(x) - (u_{0, 1}\varphi_{\epsilon_{n}})(x)|)dx \geq a_{1}, \quad \forall n \in \mathbb{N},
$$
which is absurd, because $w_{n} - u_{0, 1}\varphi_{\epsilon_{n}} \to 0$ in $L^{\Phi}(B_{\rho_1}(0))$.

Next, repeating the above arguments, we also have that $(w_{n, 2})$ given by $w_{n, 2}(x) = u_{n}^{1}(x + y_{n, 2})$ is bounded in $W^{1, \Phi}(\mathbb{R}^{N})$, and so, there exists a solution $u_{0, 2} \in W^{1, \Phi}(\mathbb{R}^{N})$ of $(P^{2})$ such that
\[
w_{n, 2}(x)\to u_{0, 2}(x)\,\, \mbox{and} \quad  \nabla w_{n, 2}(x)\to \nabla u_{0, 2}(x)\,\,\,\mbox{a. e. in } \,\,\mathbb{R}^{N}.
\]
Setting $u_{n}^{2}(x) = u_{n}^{1}(x) - (u_{0, 2}\varphi_{\epsilon_{n}})(x - y_{n, 2})$ and  arguing as above, it follows that
\[
J_{\epsilon_{n}}(u_{n}^{2})\rightarrow c - J_{0, 1}(u_{0, 1}) -J_{0, 2}(u_{0, 2}) \quad \mbox{and}\quad \Vert J^{'}_{\epsilon_{n}}(u_{n}^{2})\Vert_{\epsilon_{n}}^{*}\rightarrow 0,
\]
showing that $(u_{n}^{2})$ is a $(PS)_{ c - J_{0, 1}(u_{0, 1}) -J_{0, 2}(u_{0, 2})}^{*}$ sequence.
Continuing with this argument, we find a sequence $(u_{n}^{s})$ given by
\[
u_{n}^{s}(x) = u_{n}^{s-1}(x) - (u_{0, s}\varphi_{\epsilon_{n}})(x - y_{n, s})
\]
with 
\[
J_{\epsilon_{n}}(u_{n}^{s})\rightarrow c - \sum_{i=1}^{s}J_{0, i}(u_{0, i}) \quad \mbox{and}\quad \Vert J^{'}_{\epsilon_{n}}(u_{n}^{s})\Vert_{\epsilon_{n}}^{*}\rightarrow 0.
\]
Finally, as in \cite[Proposition 2.2]{G}, there exists $p \in \mathbb{N}$ such that
\[
J_{\epsilon_{n}}(u_{n}^{p})\rightarrow 0 \quad \mbox{and}\quad \Vert J^{'}_{\epsilon_{n}}(u_{n}^{p})\Vert_{\epsilon_{n}}^{*}\rightarrow 0.
\]
This implies that
\[
\Vert u_{n}^{p} \Vert_{\epsilon_{n}} \to 0 \quad \quad  \mbox{and}\quad \quad c = \sum_{i=1}^{p}J_{0, i}(u_{0, i}),
\]
finishing the proof.
\end{pf}
\section{Existence of a special solution for $(\tilde{P}_\epsilon)$}
Our goal is looking for a special critical point of $J_{\epsilon}$ for $\epsilon$ small enough, which will help us to prove the existence of multi-peak solutions for $(P_\epsilon)$.

In what follows, for each $i \in \Gamma$, $\widetilde{\Omega}_{\epsilon, i}$ denote mutually disjoint open sets compactly containing $\Omega_{\epsilon, i}$. Hereafter, let us denote by $E_{i}: W^{1, \Phi}(\mathbb{R}^{N}) \rightarrow \mathbb{R}$ and $E_{\epsilon, i}: \widetilde{X}_{\epsilon, i} \rightarrow \mathbb{R}$ the following functionals
\[
E_{i}(u) = \int_{\mathbb{R}^{N}}\Phi(|\nabla u|)dx + \int_{\mathbb{R}^{N}}\alpha_{i}\Phi(|u|)dx - \int_{\mathbb{R}^{N}}F(u)dx
\]
and
\[
\widetilde{E}_{\epsilon, i}(u) = \int_{\widetilde{\Omega}_{\epsilon, i}}\Phi(|\nabla u|)dx + \int_{\widetilde{\Omega}_{\epsilon, i}}V(\epsilon x)\Phi(|u|)dx - \int_{\widetilde{\Omega}_{\epsilon, i}}G(\epsilon x, u)dx,
\]
where $\widetilde{X}_{\epsilon, i}$ denotes the space of $W^{1, \Phi}(\widetilde{\Omega}_{\epsilon, i})$ 
endowed with the norm
\[
\Vert u\Vert_{\widetilde{X}_{\epsilon, i}} = \Vert \nabla u\Vert_{\Phi, \widetilde{\Omega}_{\epsilon, i}} + \Vert u\Vert_{\Phi, V_{\epsilon}, \widetilde{\Omega}_{\epsilon, i}}
\]
where
$$
\Vert \nabla u \Vert_{\Phi, \widetilde{\Omega}_{\epsilon, i}}:= \inf\Big\{ \lambda > 0;  \int_{\widetilde{\Omega}_{\epsilon, i}}\Phi\Big(\frac{\vert \nabla u\vert}{\lambda}\Big)dx \leq 1\Big\}
$$
and
\[
\Vert u\Vert_{\Phi, V_{\epsilon}, \widetilde{\Omega}_{\epsilon, i}}:= \inf\Big\{ \lambda > 0;  \int_{\widetilde{\Omega}_{\epsilon, i}}V(\epsilon x)\Phi\Big(\frac{\vert u\vert}{\lambda}\Big)dx \leq 1\Big\}.
\]

The same type of arguments found in \cite{AS} and \cite{AS2} guarantee the existence of functions $w_i \in W^{1, \Phi}(\mathbb{R}^{N})$ and $w_{\epsilon, i} \in \widetilde{X}_{\epsilon, i}$ with
\[
E_{i}(w_{i}) = \mu_{i}, \widetilde{E}_{\epsilon, i}(w_{\epsilon, i}) = \widetilde{\mu}_{\epsilon, i} \quad \mbox{and}  \quad E_{i}^{'}(w_{i})=\widetilde{E}_{\epsilon, i}^{'}(w_{\epsilon, i}) = 0,
\]
where
\[
\mu_{i} = \inf_{u \in W^{1, \Phi}(\mathbb{R}^{N})\backslash \{0\}}\sup_{t\geq 0}E_{i}(tu)= \inf_{\alpha \in \Gamma_{i}}\sup_{t \in [0, 1]}E_{i}(\alpha(t)),
\]
\[
\widetilde{\mu}_{\epsilon, i} = \inf_{u \in \widetilde{X}_{\epsilon, i} \{0\}}\sup_{t\geq 0}E_{i}(tu)= \inf_{\alpha \in \widetilde{\Gamma}_{\epsilon, i}}\sup_{t \in [0, 1]}E_{i}(\alpha(t)),
\]
\[
\Gamma_{i}=\big\{ \alpha \in C([0,1], W^{1, \Phi}(\mathbb{R}^{N}))\, :\, \alpha(0)=0,\, E_{i}(\alpha(1))<0\big\}
\]
and
\[
\widetilde{\Gamma}_{\epsilon, i}=\big\{ \alpha \in C([0,1], \widetilde{X}_{\epsilon, i})\, :\, \alpha(0)=0,\, \widetilde{E}_{\epsilon, i}(\alpha(1))<0 \big\}.
\]
\subsection{Some results about the minimax levels}

The main goal this subsection is to show an important limit involving the numbers $\mu_i$ and $\tilde{\mu}_{\epsilon,i}$. 

\begin{lem} \label{Nehari3}
	For each $i \in \Gamma$, there exist $\sigma_0, \sigma_1>0$, independents of $\epsilon$, such that
	\begin{eqnarray*}
		\Vert u\Vert_{\widetilde{X}_{\epsilon, i}} > \sigma_0\quad \mbox{and}\quad \widetilde{E}_{\epsilon, i}(u)> \sigma_{1}, \quad \forall u \in \widetilde{\mathcal{N}}_{\epsilon, i}
	\end{eqnarray*}
where
$$
\widetilde{\mathcal{N}}_{\epsilon, i} = \big\{ u \in \widetilde{X}_{\epsilon, i}\backslash \{0 \} \ : \ \widetilde{E}^{'}_{\epsilon, i}(u)u = 0 \big\}.
$$

\end{lem}
\begin{pf}
Note that, for any $u \in \widetilde{\mathcal{N}}_{\epsilon, i}$, the conditions \ref{H2} and \ref{g23} imply that 
\[
c_{1}\Big[\int_{\widetilde{\Omega}_{\epsilon, i}}\Phi(\vert \nabla u\vert)dx + \int_{\widetilde{\Omega}_{\epsilon, i}}V(\epsilon x)\Phi(\vert u\vert)dx\Big]\leq c_{2}\int_{\widetilde{\Omega}_{\epsilon, i}}\Phi_{*}(\vert u \vert)dx
\]
for some positive constants $c_{1}, c_{2}>0$. The last inequality together with \cite[Lemmas 2.3 and 2.5]{AGJ} leads to
\[
c_{1}\big(\xi_{0}(\Vert \nabla u\Vert_{\Phi, \widetilde{\Omega}_{\epsilon, i}})+ \xi_{0}(\Vert u\Vert_{\Phi, V_{\epsilon}, \widetilde{\Omega}_{\epsilon, i}})\big)\leq c_{2}\xi_{3}(\Vert u\Vert_{\Phi_{*}, \widetilde{\Omega}_{\epsilon, i}})
\]
where $\xi_0(t)=\min\{t^{l},t^{m}\}$ and $\xi_3(t)=\min\{t^{l^*},t^{m^*}\}$. Then, by Proposition \ref{TEmb} ( see Appendix ), there exists a positive constant $M^{*}$, independent of $\epsilon$, such that
\[
c_{1}\big(\xi_{0}(\Vert \nabla u\Vert_{\Phi, \widetilde{\Omega}_{\epsilon, i}})+ \xi_{0}(\Vert u\Vert_{\Phi, V_{\epsilon}, \widetilde{\Omega}_{\epsilon, i}})\big)\leq c_{2}M^{*}\xi_{3}(\Vert u\Vert_{\widetilde{X}_{\epsilon, i}}).
\]
From this, there is $\sigma_0>0$ satisfying
\[
\Vert u\Vert_{\widetilde{X}_{\epsilon, i}} > \sigma_{0}, \quad \forall u \in \widetilde{\mathcal{N}}_{\epsilon, i}.
\]
On the other hand, for any $u \in \widetilde{\mathcal{N}}_{\epsilon, i}$, the conditions \ref{H2} and \ref{g13}-\ref{g23} give
$$
\widetilde{E}_{\epsilon, i}(u) \geq  C\big( \xi_{0}(\Vert \nabla u\Vert_{\Phi, \widetilde{\Omega}_{\epsilon, i}}) + \xi_{0}(\Vert u\Vert_{\Phi, V_{\epsilon}, \widetilde{\Omega}_{\epsilon, i}})\big)
$$
for some  positive constant $C$. Therefore,
\[
\widetilde{E}_{\epsilon, i}(u)\geq \sigma_{1},
\]
for some $\sigma_1>0$. This proves the lemma. 
\end{pf}

\vspace{0.5 cm}

Our next result studies the behavior of the minimax level $\tilde{\mu}_{\epsilon,i}$ when $\epsilon$ goes to zero.
\begin{lem}\label{Lmu1}
For each $i \in \Gamma$, the following limit holds
\[
\widetilde{\mu}_{\epsilon, i}\rightarrow \mu_{i}\quad \mbox{as}\,\, \epsilon\rightarrow 0.
\]
\end{lem}
\begin{pf}
To begin with, let us prove that
\begin{eqnarray}\label{AF3}
\widetilde{\mu}_{\epsilon, i}\leq \mu_{i} + o(\epsilon).
\end{eqnarray}
In what follows, let $w_{i} \in W^{1, \Phi}(\mathbb{R}^{N})$ such that $E_{i}(w_{i})= \mu_{i}$ and $E^{'}_{i}(w_{i}) = 0$. For $\delta>0$ enough small, we fix $\vartheta \in C^{\infty}_{0}\big([0, +\infty), [0, 1]\big)$ with $\vartheta(s)= 1$ if $s \in [0, \frac{\delta}{2}]$ and $\vartheta(s)= 0$ if $s \in [\delta, +\infty)$. Using the function $\vartheta$, we define
\[
w_{\epsilon, i}(x)=\vartheta(|\epsilon x - x_{i}|)w_{i}(\frac{\epsilon x - x_{i}}{\epsilon}),
\]
where $V(x_{i}) = \displaystyle\min_{y \in \overline{\Omega}_{i}}V(y)$. As supp$(w_{\epsilon, i})\subset B_{\delta}(\frac{x_{i}}{\epsilon})$, we derive that $w_{\epsilon, i} \in \widetilde{X}_{\epsilon, i}$. Furthermore, there exists $t_{\epsilon, i}>0$ such that $\Psi_{\epsilon_{n}, i}:=t_{\epsilon, i}w_{\epsilon, i} \in \widetilde{\mathcal{N}}_{\epsilon, i}$ and
\begin{eqnarray}\label{AF4}
\widetilde{\mu}_{\epsilon, i} \leq \max_{t\geq 0}\widetilde{E}_{\epsilon, i}(tw_{\epsilon, i})= \widetilde{E}_{\epsilon, i}(t_{\epsilon, i}w_{\epsilon, i}).
\end{eqnarray}
Using Lebesgue's Theorem, it is possible to prove that
$$
\lim_{\epsilon \rightarrow 0}\widetilde{E}_{\epsilon, i}(t_{\epsilon_{n}, i}w_{\epsilon_{n}, i}) = E_{i}(w_{i})=\mu_i.
$$
Consequently,
\begin{eqnarray}\label{AF5}
\limsup_{\epsilon \to 0}\widetilde{\mu}_{\epsilon, i} \leq \mu_i.
\end{eqnarray}

Now, we will prove the inequality below
\begin{eqnarray}\label{AF6}
\mu_{i}\leq \liminf_{\epsilon \to 0} \widetilde{\mu}_{\epsilon, i}.
\end{eqnarray}
Let $\epsilon_{n}\in (0, +\infty)$ with $\epsilon_{n}\rightarrow 0$ and $v_{\epsilon_{n}, i} \in \widetilde{X}_{\epsilon_{n}, i}$ be a solution of the following problem
\begin{align}
\left\{
\begin{array}
[c]{rcl}%
- \Delta_{\Phi}u + V(\epsilon_{n} x)\phi(\vert u\vert)u & = & g(\epsilon x, u)~  \mbox{in}~ \widetilde{\Omega}_{\epsilon_{n}, i},\\
\displaystyle\frac{\partial u}{\partial \nu} = 0, ~  \mbox{on}~  \partial \widetilde{\Omega}_{\epsilon_{n}, i}. &  &
\end{array}
\right. \tag{$ P_{\epsilon, i}$}\label{PA33}%
\end{align}
By Lemma \ref{Nehari3}, there exists $\sigma_{0}>0$, independent of $n$, such that
\[
\Vert v_{\epsilon_{n}, i}\Vert_{X_{\epsilon_{n}, i}}\geq \sigma_{0}, \,\,\,\mbox{for all}\,\, n \in \mathbb{N}.
\]
Using the last inequality together with the Proposition \ref{LEA3} (see Appendix), there exist $(y_{n, i})\subset \mathbb{R}^{N}$, $\varrho>0$ and $a>0$ such that
\begin{eqnarray}\label{DS3}
\lim_{n\rightarrow +\infty}\int_{B_{\varrho}(y_{n, i})\cap \Omega_{\epsilon_n, i}}\Phi(|v_{\epsilon_{n}, i}|)dx \geq a.
\end{eqnarray}
Moreover, by \eqref{DS3}, increasing $\varrho$ if necessary, we may assume that $(y_{n, i}) \subset \Omega_{\epsilon_n, i}$ with dist$(y_{n,i}, \partial\widetilde{\Omega}_{\epsilon_n, i}) \to +\infty$. Hence, $\epsilon_{n}y_{n, i} \to \overline{x}_{i}\in \overline{\Omega}_{\epsilon_n, i}$ and given $\rho > \varrho $, we have $B_{2\rho}(y_{n, i})\subset \widetilde{\Omega}_{\epsilon_n,i}$ for $n$ sufficiently large. Setting
\[
w_{n, i, \rho}(x)= \psi\Big(\frac{|x|}{\rho}\Big)v_{\epsilon_{n}, i}(x + y_{n, i}), \quad \forall x \in \widetilde{\Omega}_{\epsilon_{n}, i} - y_{n, i}
\]
where $\psi \in C^{\infty}(\mathbb{R})$ is such that $\psi = 1$ on $[0, 1]$, $\psi = 0$ on $(2, +\infty)$, $0\leq \psi \leq 1$ and $\psi^{'} \in L^{\infty}(\mathbb{R})$, we find
\[
\int_{B_{\varrho}(0)}\Phi(|w_{n, i, \rho}|)dx = \int_{B_{\varrho}(y_{n, i})}\Phi(|v_{\epsilon_{n}, i}|)dx\geq a>0.
\]
Once supp$(w_{n, i, \rho})\subset B_{2\rho}(0)$, we conclude  that $w_{n, i, \rho} \in W^{1, \Phi}(\mathbb{R}^{N})$.  The fact that $v_{\epsilon_{n}, i}$  is a solution of \eqref{PA33} together with \eqref{AF3} yields there exists $C>0$, independent of $\rho$, such that $\Vert w_{n, i, \rho}\Vert\leq C$. Hence, there exists $w_{\rho}^{i}\in W^{1, \Phi}(\mathbb{R}^{N})$ such that
\[
w_{n, i, \rho}\rightharpoonup w_{\rho}^{i}\,\,\,\mbox{in}\,\,W^{1, \Phi}(\mathbb{R}^{N}).
\]
Then, $w_{n, i, \rho}\to w_{\rho}^{i}$ in $L^{\Phi}_{loc}(\mathbb{R}^{N})$ and 
\begin{eqnarray}
\int_{B_{\varrho}(0)}\Phi(|w_{\rho}^{i}|)dx\geq a>0.
\end{eqnarray}
Since $(\Vert w_{\rho}^{i}\Vert)$ is bounded in $\mathbb{R}$, there exists $w\in W^{1, \Phi}(\mathbb{R}^{N})$ such that
\[
w_{\rho}^{i}\rightharpoonup w^{i}\,\,\,\mbox{in}\,\,W^{1, \Phi}(\mathbb{R}^{N}).
\]
Thus, $w_{\rho}^{i}\to w^{i}$ in $L^{\Phi}_{loc}(\mathbb{R}^{N})$ and 
\begin{eqnarray}
\int_{B_{\varrho}(0)}\Phi(|w^{i}|)dx\geq a>0.
\end{eqnarray}
Moreover, by a direct computation, $w^{i}$ is a solution of problem $(P^{i})$, that is,
\[
\int_{\mathbb{R}^{N}}\Big[\phi(|\nabla w^{i}|)\nabla w^{i}\nabla \zeta dx + V(\overline{x}_{i}) \phi(|w^{i}|)w^{i}\zeta\Big] dx = \int_{\mathbb{R}^{N}}g(\overline{x}_{i}, w^{i})w^{i}\zeta dx,
\]
for all $\zeta \in W^{1,\Phi}(\mathbb{R}^N)$. Fixing $\tau>\rho$, we know that $B_{\tau}(y_{n, i})\subset \widetilde{\Omega}_{\epsilon_{n}, i}$ for $n$ large enough. Hence,
\begin{eqnarray*}
\widetilde{\mu}_{\epsilon_{n}, i} &=& \widetilde{E}_{\epsilon_{n}, i}(v_{\epsilon_{n}, i}) -\frac{1}{\theta}\widetilde{E}^{'}_{\epsilon_{n}, i}(v_{\epsilon_{n}, i})v_{\epsilon_{n}, i}\\
&\geq&\int_{B_{\tau}(0)}\Big[h(|\nabla w_{n, i, \rho}|) + V(\epsilon_{n}x + \epsilon_{n}y_{n, i})h(| w_{n, i, \rho}|)\Big]dx\\
&& + \int_{B_{\tau}(0)}\Big[\frac{1}{\theta}g(\epsilon_{n}x, w_{n, i, \rho})w_{n, i, \rho} - G(\epsilon_{n}x,w_{n, i, \rho})\Big]dx
\end{eqnarray*}
where $h(t)= \Phi(t)-\frac{1}{\theta}\phi(t)t^{2}$. Applying the Fatou's lemma in $n$, and after taking the limit of $\rho \to +\infty$, we derive that
\begin{eqnarray*}
\liminf_{n\to +\infty}\widetilde{\mu}_{\epsilon_{n}, i} &\geq&\int_{\mathbb{R}^{N}}\Big[h(|\nabla w^{i}|) + V(\overline{x}_{i})h(| w^{i}|)\Big]dx + \int_{\mathbb{R}^{N}}\Big[\frac{1}{\theta}g(\overline{x}_{i}, w^{i})w^{i} - G(\overline{x}_{i},w^{i})\Big]dx\\
&=& J_{0, i}(w^{i}) -\frac{1}{\theta}J_{0, i}^{'}(w^{i})w^{i} = J_{0, i}(w^{i}) = J_{0, i}(w^{i}) = \mu_{V(\overline{x}_{i})}\geq \mu_{i},
\end{eqnarray*}
showing \eqref{AF6}. By \eqref{AF3} and \eqref{AF6}, 
\[
\widetilde{\mu}_{\epsilon, i}\rightarrow \mu_{i}\quad \mbox{as}\,\, \epsilon\rightarrow 0,
\]
which proves the lemma.
\end{pf}
\subsection{Critical points for $J_{\epsilon}$}

In the sequel, we fix $\Gamma \subset \big\{1, ..., \kappa\big\}$ and for each $i \in \Gamma$, we choose $\rho_{i}>1$ such that $E_{i}(\rho_{i}^{-1}w_{i}), E_{i}(\rho_{i} w_{i})<\mu_{i}$. Setting $\rho = \displaystyle\max_{i \in \Gamma}\rho_{i}$, we have
\begin{eqnarray}\label{Omu}
E_{i}(\rho^{-1}w_{i}), E_{i}(\rho w_{i})<\mu_{i}\,\,\,\mbox{for all}\,\, i \in \Gamma
\end{eqnarray}
and
\[
\mu_{i}=\max_{t \in [\rho^{-2}, 1]}E_{i}(t\rho w_{i})\,\,\,\mbox{for all}\,\, i \in \Gamma.
\]

Moreover, without loss of generality, we will consider $\Gamma=\{1, ..., \lambda\}$ for some $\lambda \in \big\{1, ..., \kappa\big\}$ and define $\widetilde{H}_{\epsilon} : [\rho^{-2}, 1]^{\lambda} \to X_{\epsilon}$ by
\begin{eqnarray}\label{FunH}
\widetilde{H}_{\epsilon}(\overrightarrow{\theta})(z)= \sum_{i=1}^{\lambda}\theta_{i}\rho(w_{i}\varphi)\big(z-\frac{x_{i}}{\epsilon}\big)
\end{eqnarray}
for all $\overrightarrow{\theta}= (\theta_{1}, ...,\theta_{\lambda}) \in [\rho^{-2}, 1]^{\lambda}$, where $x_{i} \in \Upsilon_{i}= \{x \in \Omega_{i}\,\, :\,\, V(x_{i})=\alpha_{i}\}$. Moreover, we set
\begin{eqnarray*}
&&\mathcal{U}_{\epsilon}= \big\{ H \in C\big([\rho^{-2}, 1]^{\lambda}, X_{\epsilon}\big); \,\, H = \widetilde{H}_{\epsilon}\,\,\mbox{on}\,\, \partial ([\rho^{-1}, 1]^{\lambda})\big),\,\,\\
&&\quad \quad\quad \quad\quad \quad H(\overrightarrow{\theta})\mid_{\Omega_{\epsilon, i}}\neq 0\,\,\forall i \in \Gamma\,\,\mbox{and}\,\, \forall\overrightarrow{\theta} \in [\rho^{-1}, 1]^{\lambda}\big\}.
\end{eqnarray*}
Since $\mbox{supp}\Big(w_{i}\varphi\big(z-\frac{x_{i}}{\epsilon}\big)\Big)\subset \Omega_{\epsilon, i}$, it follows that $\widetilde{H}_{\epsilon} \in \mathcal{U}_{\epsilon}$. Therefore, we can define the number
\[
\mathcal{S}_{\epsilon}=\inf_{H \in \mathcal{U}_{\epsilon}}\max_{\overrightarrow{\theta} \in [\rho^{-2}, 1]^{\lambda}}J_{\epsilon}(H(\overrightarrow{\theta})).
\]
\begin{lem}\label{Lmu2}
For $\epsilon$ small enough, the following property holds: If $H \in \mathcal{U}_{\epsilon}$, then there exists $\overrightarrow{\theta_{*}} \in [\rho^{-2}, 1]^{\lambda}$, such that
\[
\widetilde{E}^{'}_{\epsilon, i}\big(H(\overrightarrow{\theta_{*}})\big)H(\overrightarrow{\theta_{*}}) = 0,\,\,\,\mbox{for all}\,\, i \in \Gamma.
\]
In particular, $\widetilde{E}_{\epsilon, i}\big(H(\overrightarrow{\theta_{*}})\big) \geq \widetilde{\mu}_{\epsilon, i}, \,\,\, i = 1, ..., \lambda$.
\end{lem}
\begin{pf}
Given $H \in \mathcal{U}_{\epsilon}$, consider $\overline{H}: [\rho^{-2}, 1]^{\lambda} \to \mathbb{R}^{\lambda}$ such that
\[
\overline{H}(\overrightarrow{\theta}) = \Big(\widetilde{E}^{'}_{\epsilon, 1}\big(H(\overrightarrow{\theta})\big)H(\overrightarrow{\theta}), ...., \widetilde{E}^{'}_{\epsilon, \lambda}\big(H(\overrightarrow{\theta})\big)H(\overrightarrow{\theta})\Big),\,\,\,\mbox{where}\,\, \overrightarrow{\theta}= (\theta_{1}, ...,\theta_{\lambda}).
\]
For $\overrightarrow{\theta} \in \partial \big([\rho^{-1}, 1]^{\lambda}\big)$, it holds 
$$
\overline{H}(\overrightarrow{\theta}) = \Big(\widetilde{E}^{'}_{\epsilon, 1}\big(\widetilde{H}_{\epsilon}(\overrightarrow{\theta})\big)\widetilde{H}_{\epsilon}(\overrightarrow{\theta}), ...., \widetilde{E}^{'}_{\epsilon, \lambda}\big(\widetilde{H}_{\epsilon}(\overrightarrow{\theta}))\big)\widetilde{H}_{\epsilon}(\overrightarrow{\theta})\Big).
$$ 
From this, we observe that there is no $\overrightarrow{\theta} \in \partial \big([\rho^{-2}, 1]^{\lambda}\big)$ with $\overline{H}(\overrightarrow{\theta}) = 0$. In fact, for all $i \in \Gamma$
\[
\widetilde{E}^{'}_{\epsilon, i}(\widetilde{H}_{\epsilon}(\overrightarrow{\theta}))\widetilde{H}_{\epsilon}(\overrightarrow{\theta}) = E^{'}_{i}(\theta_{i}\rho w_{i})\theta_{i}\rho w_{i} + o_{\epsilon}(1)\,\,\,\mbox{uniformly in}\,\,\overrightarrow{\theta} \in [\rho^{-2}, 1]^{\lambda}.
\]
Thereby, if $\overrightarrow{\theta} \in \partial \big([\rho^{-2}, 1]^{\lambda}\big)$, then $\theta_{i_{0}} = 1$  or $\theta_{i_{0}} = \rho^{-2}$ for some $i_{0} \in \Gamma$. Consequently,
\[
0=\widetilde{E}^{'}_{\epsilon, i_{0}}(\overline{H}(\overrightarrow{\theta}))\overline{H}(\overrightarrow{\theta}) = E^{'}_{i_{0}}(\rho w_{i_{0}})\rho w_{i_{0}} + o_{\epsilon}(1)
\]
or
\[
0=\widetilde{E}^{'}_{\epsilon, i_{0}}(\overline{H}(\overrightarrow{\theta}))\overline{H}(\overrightarrow{\theta}) = E^{'}_{i_{0}}(\rho^{-2} w_{i_{0}})\rho^{-2} w_{i_{0}} + o_{\epsilon}(1).
\]
Therefore, if $\widetilde{E}^{'}_{\epsilon, i_{0}}(\overline{H}(\overrightarrow{\theta}))\overline{H}(\overrightarrow{\theta})=0$, the limit of $\epsilon \to 0$ gives
\[
E^{'}_{i_{0}}(\rho w_{i_{0}})\rho w_{i_{0}}=0 \quad \mbox{or} \quad E^{'}_{i_{0}}(\rho^{-2} w_{i_{0}})\rho^{-2} w_{i_{0}}=0
\]
from where it follows that
\[
E_{i_{0}}(\rho w_{i}) \geq  \mu_{i}\,\,\,\mbox{or}\,\,E_{i_{0}}(\rho^{-2} w_{i})\geq \mu_{i},
\]
Thereby, there exists $\overrightarrow{\theta_{*}} \in (\delta^{-1}, 1)^{\lambda}$ satisfying
\[
\widetilde{E}^{'}_{\epsilon, i}\big(H(\overrightarrow{\theta_{*}})\big)H(\overrightarrow{\theta_{*}}) = 0,\,\,\,\mbox{for all}\,\, i \in \Gamma.
\]
\end{pf}

The next result establishes an important relation between $\mathcal{S}_{\epsilon}$ and the levels $\mu_{i}$. In what follows, we consider $\mathcal{D}_{\Gamma} = \displaystyle\sum_{i=1}^{\lambda}\mu_{i} $. By using the same ideas found in \cite{Alves2011}, it is possible to prove the following results 
\begin{prop}\label{Mw0}
The following limit holds
\[
\lim_{\epsilon \to 0}\mathcal{S}_{\epsilon} = \mathcal{D}_{\Gamma}.
\]
\end{prop}
\begin{corol}
For each $\alpha>0$, there exists $\epsilon_{0} = \epsilon_{0}(\alpha)$ such that
\begin{eqnarray*}
\sup_{\overrightarrow{\theta} \in [\delta^{-1}, 1]^{\lambda}}J_{\epsilon}\big(\widetilde{H}_{\epsilon}(\overrightarrow{\theta})\big)\leq \mathcal{D}_{\Gamma} + \frac{\alpha}{2}\,\,\, \forall \epsilon \in (0, \epsilon_{0}).
\end{eqnarray*}
\end{corol}

Next, we will introduce some notations. Firstly, we fix the set
\[
\mathcal{Z}_{\epsilon, i} = \big\{u \in \widetilde{X}_{\epsilon, i}\,\, :\,\,\Vert u\Vert_{\widetilde{X}_{\epsilon, i}} \leq \frac{\sigma_{0}}{2}\big\}
\]
where $\sigma_{0} >0$ is a constant such that
\[
\liminf_{\epsilon \to 0}\Vert \widetilde{H}_{\epsilon}(\overrightarrow{\theta})\Vert_{\widetilde{X}_{\epsilon, i}} > \sigma_0\,\,\,\mbox{uniformly in }\,\, \overrightarrow{\theta} \in [\rho^{-2}, 1]^{\lambda}\,\,\mbox{and}\,\, i \in \Gamma \,\, ( \mbox{See Lemma \ref{Nehari3}} ). 
\]
Hence, there exist positive constants $\tau$ and $\epsilon^{*}$ such that
\[
\mbox{dist}_{\epsilon, i}\big( \widetilde{H}_{\epsilon}(\overrightarrow{\theta}), \mathcal{Z}_{\epsilon, i}\big) > \tau\,\,\,\mbox{for all}\,\, \overrightarrow{\theta} \in [\delta^{-2}, 1]^{\lambda}, \, i \in \Gamma\,\,\mbox{and}\,\, \epsilon \in (0, \epsilon^{*}),
\]
where dist$_{\epsilon, i}(A, B)$ denotes the distance between sets $A$ and $B$ of $\widetilde{X}_{\epsilon, i}$. Moreover, we define
\[
\Theta = \big\{u \in X_{\epsilon}\, :\, \mbox{dist}_{\epsilon, i}(u, \mathcal{Z}_{\epsilon, i} )\geq \tau\,\, \mbox{for all}\,\, i \in \Gamma\big\}
\]
and for any $c, \mu>0$ and $0< \delta<\frac{\tau}{2}$, we consider the sets
\[
J_{\epsilon}^{c} =\big\{u \in X_{\epsilon}\,\, :\,\, J_{\epsilon}(u)\leq c\big\}\,\,\,\mbox{and}\,\, \mathcal{Q}_{\epsilon, \mu} = \big\{ u \in \Theta_{2\delta}\,\, :\,\, \vert J_{\epsilon}(u) - \mathcal{S}_{\epsilon}\vert\leq \mu\big\},
\]
where $\Theta_{s}$, for $s>0$, denotes the set
\[
\Theta_{s} = \big\{u \in X_{\epsilon}\,\, :\,\, \mbox{dist}(u, \Theta) \leq s\big\}.
\]
Observe that for each $\mu>0$, there exists $\epsilon_{1}=\epsilon_{1}(\mu)>0$ such that the function $U_{ \epsilon}$ given by
\[
U_{\epsilon}(z)= \sum_{i=1}^{\lambda}\Big(w_{i}\varphi\big(z-\frac{x_{i}}{\epsilon}\big)\Big)
\]
verifies
\[
U_{\epsilon}\in \Theta_{s}\,\,\,\mbox{for all}\,\, \epsilon \in (0, \epsilon_{1})\,\,\,\mbox{and}\,\,J_{\epsilon}(U_{\epsilon}) = \mathcal{D}_{\Gamma} + o_{\epsilon}(1).
\]
As $\mathcal{S}_{\epsilon} = \mathcal{D}_{\Gamma} + o_{\epsilon}(1)$, we have
\[
J_{\epsilon}(U_{\epsilon}) = \mathcal{S}_{\epsilon} + o_{\epsilon}(1),
\]
showing that $\mathcal{Q}_{\epsilon, \mu}\neq\emptyset$.

Next, let us consider $M$ large enough, independent of $\epsilon$, satisfying
\begin{eqnarray}\label{dpa1}
\Vert \widetilde{H}_{\epsilon}(\overrightarrow{\theta})\Vert_{\epsilon}\leq \frac{M}{2}\,\,\,\mbox{for all}\,\, \overrightarrow{\theta} \in [\rho^{-2}, 1]^{2\lambda}.
\end{eqnarray}
For each $s>0$, we denote by $ \overline{B}_{s} = \big\{ u \in X_{\epsilon}\, ;\, \Vert u\Vert_{\epsilon}\leq s\big\} $ and define the number
\[
\mu_{*} =\min\Big\{\frac{\mu_{i}}{4}, \frac{M}{4}, \frac{\delta}{4}; \,\, i \in \Gamma \Big\}.
\]

The result below establishes the existence of a special critical point for functional $J_{\epsilon}$, which will be used later on. However, we will omit its proof because it follows by using the same approach explored in \cite{Alves2011}.
\begin{prop}\label{PCJ}
For each $\mu \in (0, \mu_{*})$, there exists $\epsilon_{\mu}>0$ such that $J_{\epsilon}$ has a critical point $v_{\epsilon} \in \mathcal{Q}_{\epsilon, \mu}\cap \overline{B}_{M+1}\cap J_{\epsilon}^{\mathcal{D}_{\Gamma}}$ for all $\epsilon \in (0, \epsilon_{\mu})$.
\end{prop}

\section{The existence of multi-peak positive solutions}
In this section, we will show existence of $\lambda$-peak solution for \eqref{PP3}. For this purpose, we need of the following technical lemma

\begin{lem}\label{LemF3}
There exist $\overline{\epsilon}, \overline{\mu}$, such that the solution $v_{\epsilon}$ obtained in Proposition \ref{PCJ} satisfies
\[
\max_{z \in \partial \Omega_{\epsilon}}v_{\epsilon}(z) < a\,\,\,\mbox{for all}\,\, \mu \in (0, \overline{\mu})\,\,\mbox{and}\,\, \epsilon \in (0, \overline{\epsilon}).
\]
\end{lem}
\begin{pf}
Assume by contradiction, that there exist $\epsilon_{n}, \mu_{n} \to 0$ such that
\[
v_{n} :=v_{\epsilon_{n}} \in \mathcal{Q}_{\epsilon_{n}, \mu_{n}}\,\,\,\mbox{and}\,\,\,\max_{z \in \partial \Omega_{\epsilon_{n}}}v_{n}(z)\geq a\,\,\,\mbox{for all}\,\, n \in \mathbb{N}.
\]
Since $v_{n} \in \mathcal{Q}_{\epsilon_{n}, \mu_{n} }$, we know that 
\begin{equation} \label{Li23'}
J^{'}_{\epsilon_{n}}(v_{n}) = 0, \quad |J_{\epsilon_{n}}(v_{n}) - \mathcal{S}_{\epsilon_{n}} |\to 0 \,\,\,\mbox{and}\,\,\,\mbox{dist}(v_{n}, \Theta)\leq 2\delta.
\end{equation}
Applying the Proposition \ref{PropI3}, there exist a nonnegative integer $p$, sequences of points $(y_{n, i})\subset \mathbb{R}^{N}$, points $x_{i} \in \overline{\Omega}$, $i = 1, ..., \lambda$ and functions $u_{0, i}$ verifying 
\begin{equation}\label{Li23}
\Big\Vert v_{n}(\cdot)-\sum_{i=1}^{p}u_{0, i}(\cdot - y_{n, i})\varphi_{\epsilon_{n}}(\cdot - y_{n, i})\Big\Vert_{\epsilon_{n}} \rightarrow 0\,\,\,\mbox{as}\,\,n \rightarrow +\infty
\end{equation}
and
\begin{equation}\label{Li13}
\epsilon_{n}y_{n, i}\rightarrow x_{i}\,\,\, \mbox{for}\,\, i = 1, ..., p.
\end{equation}
From \eqref{Li23'}, \eqref{Li23} and \eqref{Li13},  $p = \lambda$ and $x_{i} \in \Upsilon_{i}$ for all $i = 1, ..., \lambda$. 

In what follows, we fix $(z_{n}) \subset \partial \Omega_{\epsilon_{n}}$ such that
\[
v_{n}(z_{n}) = \max_{z \in \partial \Omega_{\epsilon}}v_{\epsilon_{n}}(z)
\]
and the function $w_{n}(x) = v_{n}(x+ z_{n})$. Then,  
\begin{eqnarray*}
\Big\Vert w_{n}(\cdot)-\sum_{i=1}^{p}u_{0, i}(\cdot + z_{n} - y_{n, i})\varphi_{\epsilon_{n}}(\cdot + z_{n} - y_{n, i})\Big\Vert_{W^{1, \Phi}(\mathbb{R}^{N})} \rightarrow 0\,\,\,\mbox{as}\,\,n \rightarrow +\infty
\end{eqnarray*}
On the other hand, for each $\varrho>0$,
\begin{eqnarray*}
	\Big\Vert \sum_{i=1}^{p}u_{0, i}(\cdot + z_{n} - y_{n, i})\varphi_{\epsilon_{n}}(\cdot + z_{n} - y_{n, i})\Big\Vert_{W^{1, \Phi}(B_{\varrho}(0))} \rightarrow 0\,\,\,\mbox{as}\,\,n \rightarrow +\infty.
\end{eqnarray*}
Consequently,
\begin{eqnarray}\label{ADD3}
\Vert w_{n}\Vert_{W^{1, \Phi}(B_{\varrho}(0))} \rightarrow 0\,\,\,\mbox{as}\,\,n \rightarrow +\infty.
\end{eqnarray}
Notice that $w_{n}$ is solution of problem
\[
- \Delta_{\Phi}w_{n} + V(\epsilon_{n} x + \epsilon_{n}z_{n})\phi(|w_{n}|)w_{n}  =  g(\epsilon_{n} x + \epsilon_{n}z_{n}, w_{n})\,\,\,\mbox{in}\,\, \mathbb{R}^{N},
\]
because $v_{n}$ is a solution of \eqref{PA3}. Arguing as in \cite[Lemma 3.2]{AS2}, there exists $w \in C_{loc}^{1, \alpha}(\mathbb{R}^{N})$ such that, up to a subsequence,
\[
w_{n}\to w\,\,\,\mbox{in}\,\, C_{loc}^{1, \alpha}(\mathbb{R}^{N}).
\]
Since
\[
\max_{z \in \partial \Omega_{\epsilon_{n}}}v_{n}(z)\geq a,
\]
we have that $w_{n}(0)\geq a$ for all $n \in \mathbb{N}$, and so, $w(0)\geq a$. Thereby, there exists $\varrho \in (0, 1)$ such that $w(x)\geq \frac{a}{2}$ for all $x \in B_{\varrho}(0)$. Consequently $w\neq 0$, which is a contradiction with \eqref{ADD3}, showing the lemma.
\end{pf}

\subsection{Proof of Theorem 1.1}
By Lemma \ref{LemF3}, there exist $\overline{\epsilon}, \overline{\mu}>0$, such that the solution $v_{\epsilon} \in \mathcal{Q}_{\epsilon, \mu }$ obtained in Proposition \ref{PCJ} satisfies
\[
\max_{z \in \partial \Omega_{\epsilon}}v_{\epsilon}(z) < a\,\,\,\mbox{for all}\,\, \mu \in (0, \overline{\mu})\,\,\mbox{and}\,\, \epsilon \in (0, \overline{\epsilon}).
\]
Repeating the same arguments found in \cite{AS2}, we see that
\[
v_{\epsilon}(x)\leq a \,\,\,\mbox{for all}\,\, x \in \mathbb{R}^{N}\backslash \Omega_{\epsilon}.
\]
Hence, $v_{\epsilon}$ is a solution of \eqref{pp3} for all $ \epsilon \in (0, \overline{\epsilon})$. To finish the proof, we will show that the family $(v_{\epsilon})$ is a $\lambda$-peak solution. To see why, we consider $\epsilon_{n}\to 0$ and $v_{n}= v_{\epsilon_{n}}$. Observe that $(v_{n})$ is a $(PS)_{\mathcal{D}_{\Gamma}}^{*}$ sequence verifying
\begin{eqnarray}\label{DT331}
\mbox{dist}(v_{n}, \Theta)\leq 2\delta\,\,\,\mbox{for all}\,\, n \in \mathbb{N}.
\end{eqnarray}
From Proposition \ref{PropI3}, there exist a subsequence of $(v_{n})$, still denoted by itself, a nonnegative integer $p$, sequences of points $(y_{n, i})\subset \mathbb{R}^{N}$ with $i = 1, ..., p$ such that
\begin{eqnarray} \label{yni}
\epsilon_{n}y_{n, i}\rightarrow x_{i} \in \overline{\Omega}\quad \mbox{and}\quad \vert y_{n, j} - y_{n, i}\vert \rightarrow +\infty\quad \mbox{as} \,\,\,n\rightarrow +\infty
\end{eqnarray}
with $i \in \big\{1, ..., p\big\}$ and
\begin{eqnarray}\label{fer}
\Big\Vert v_{n}(\cdot)-\sum_{i=1}^{p}u_{0, i}(\cdot - y_{n, i})\varphi_{\epsilon_{n}}(\cdot - y_{n, i})\Big\Vert_{\epsilon_{n}} \rightarrow 0\,\,\,\mbox{as}\,\,n \rightarrow +\infty
\end{eqnarray}
where $\varphi_{\epsilon}(x) = \varphi\big(x/ (-\ln\epsilon)\big)$ for $0 < \epsilon < 1$, and $\varphi$ is a cut-off function which $\varphi(z) = 1$ for $|z|\leq 1$, $\varphi(z) = 0$ for $|z|\geq 2$ and $\vert \nabla \varphi\vert \leq 2$. The function $u_{0, j} \neq 0$ is a nonnegative solution for
\[
- \Delta_{\Phi}u + V_{i}\phi(\vert u\vert)u  =  g_{0, i}(x, u)\quad \mbox{in} \quad \mathbb{R}^{N},
\]
where $V_{i} = V(x_{i})\geq V_{0}>0$ and $g_{0,i}(x, u) = \displaystyle\lim_{n \rightarrow \infty}g(\epsilon_{n}x + \epsilon_{n}y_{n, i}, u)$. Furthermore,
\begin{eqnarray}\label{DT334}
\sum_{i= 1}^{\lambda}\mu_{i} = \sum_{i = 1}^{p}J_{0, i}(u_{0, i})
\end{eqnarray}
where $J_{0, i} : W^{1, \Phi}(\mathbb{R}^{N}) \rightarrow \mathbb{R}$ denotes the functional given by
\[
J_{0, i}(u) = \int_{\mathbb{R}^{N}}\Phi(\vert\nabla u\vert)dx + V_{i}\int_{\mathbb{R}^{N}}\Phi(\vert u\vert)dx - \int_{\mathbb{R}^{N}}G_{0, i}(x, u)dx
\]
with $G_{0, i}(x, t) = \int_{0}^{t}g_{0, i}(x, s)ds$.
Arguing as in proof of Lemma \ref{LemF3} and using \eqref{DT331}-\eqref{DT334}, we infer that $p = \lambda$, $x_{i} \in \overline{\Omega}_{i}$ and 
\begin{eqnarray*}
\sum_{j = 1}^{\lambda}\mu_{j} = \sum_{j = 1}^{\lambda}J_{0, j}(u_{0, j}).
\end{eqnarray*}
The last equality yields $x_{i} \in \Upsilon_{i}$ and $V(x_{i})=\alpha_i$, because if for some $i_{0} \in {1, ..., \lambda}$, we have $x_{i_{0}} \in \partial \Omega_{i}$, the assumption $(V_{1})$ leads to $V(x_{i_{0}})>\alpha_{i}$, and so, $J_{0, i_{0}}(u_{0, i_{0}})>\mu_{i_{0}}$.

On the other hand, since $J_{0, i}(u_{0, i_{0}})\geq \mu_{i}$, for all $i=1, ..., \lambda$, we must have 
\begin{eqnarray*}
\sum_{i= 1}^{\lambda}\mu_{i} < \sum_{i = 1}^{p}J_{0, i}(u_{0, i}),
\end{eqnarray*}
which is a contradiction. Therefore, $V(x_{i}) = \alpha_{i}$ for $i=1, ..., \lambda$ and $u_{0, i}$ is a nontrivial solution of problem
\[
- \Delta_{\Phi}u + \alpha_{i}\phi(\vert u\vert)u  =  f(u)\quad \mbox{in} \quad \mathbb{R}^{N}.
\]
Now, we will show that for each $\eta > 0$, there exists $\rho>0$ such that
\begin{eqnarray}
\Vert v_{n} \Vert_{\infty, \mathbb{R}^{N}\backslash \cup_{j = 1}^{p}B_{\rho}(y_{n, i})} \leq  \eta
\end{eqnarray}
and  there exists $\delta>0$ such that
\begin{eqnarray}
\Vert v_{n} \Vert_{\infty, B_{\rho}(y_{n, j})}\geq \delta,\,\,\,\mbox{for all}\,\, j \in \Gamma.
\end{eqnarray}
To this end, we need of the following estimate:
\begin{claim}\label{fer1}
Given $\eta >0$, there exist $\rho>0 $ and $n_0 \in \mathbb{N}$ such that
\begin{eqnarray}
\Vert v_{n}\Vert_{W^{1, \Phi}(\mathbb{R}^{N}\backslash \cup_{j = 1}^{p}B_{\rho}(y_{n, i}))}\leq \eta, \quad \forall n \geq n_0.
\end{eqnarray}
\end{claim}
In fact, for each $j \in \Gamma$, there exists $\rho_{j}>0$ such that
\[
\Vert u_{0,j}\Vert_{W^{1, \Phi}(\mathbb{R}^{N}\backslash B_{\rho_{j}}(0))}< \eta.
\]
Setting $\rho = \max\{\rho_{1}, ..., \rho_{p}\}$, we have
\[
\int_{\mathbb{R}^{N}\backslash B_{\rho}(0)}\Phi(|\nabla u_{0,j}|)dx,\,\,\int_{\mathbb{R}^{N}\backslash B_{\rho}(0)}\Phi(|u_{0,j}|)dx < \eta\,\,\, \mbox{for all}\,\, j \in \Gamma.
\]
Notice that
\[
\int_{\mathbb{R}^{N}\backslash B_{\rho}(y_{n, j})}\Phi\big(|\nabla(u_{0, j}(\cdot - y_{n, i})\varphi_{\epsilon_{n}}(\cdot - y_{n, j}))|\big)dx = \int_{\mathbb{R}^{N}\backslash B_{\rho}(0)}\Phi\big(|\nabla (u_{0, j}\varphi_{\epsilon_{n}})|\big)dx .
\]
From  $\Delta_{2}$-condition, we get
\[
\int_{\mathbb{R}^{N}\backslash B_{\rho}(0)}\Phi\big(|\nabla (u_{0, j}\varphi_{\epsilon_{n}})|\big)dx \leq c_{1}\int_{\mathbb{R}^{N}\backslash B_{\rho}(0)}\Phi(|\nabla u_{0,j}|)dx + c_{2}\int_{\mathbb{R}^{N}\backslash B_{\rho}(0)}\Phi(|u_{0,j}|)dx.
\]
Thereby, given $\eta>0$, we can find $\rho$ large enough verifying
\[
\int_{\mathbb{R}^{N}\backslash B_{\rho}(y_{n, j})}\Phi\big(|\nabla(u_{0, j}(\cdot - y_{n, i})\varphi_{\epsilon_{n}}(\cdot - y_{n, j}))|\big)dx < \frac{\eta}{2}.
\]
Similarly,
\[
\int_{\mathbb{R}^{N}\backslash B_{\rho}(y_{n, j})}\Phi(|u_{0, j}(\cdot - y_{n, i})\varphi_{\epsilon_{n}}(\cdot - y_{n, j})|)dx < \frac{\eta}{2},
\]
showing that
\begin{eqnarray} \label{NOVAESTI}
\Vert u_{0, j}(\cdot - y_{n, j})\varphi_{\epsilon_{n}}(\cdot - y_{n, j})\Vert_{W^{1, \Phi}(\mathbb{R}^{N}\backslash B_{\rho}(y_{n, j}))}\leq \eta.
\end{eqnarray}
Now, the claim follows from \eqref{fer} and (\ref{NOVAESTI}).

Using the above information, we are able to prove the following estimate
\begin{claim}\label{afir989} Given $\eta>0$, there are $\rho>0$ and $n_0 \in\mathbb{N}$ such that
	\[
	|v_{n}(z)| \leq  \eta\,\,\, \mbox{for all}\,\, z \in \mathbb{R}^{N}\backslash \cup_{j = 1}^{p}B_{\rho +1}(y_{n, i}), \quad \forall n \geq n_0.
	\]
\end{claim}
Indeed, fix $R_1 \in (0, 1)$ and $x_{0} \in \mathbb{R}^{N}\backslash \cup_{j = 1}^{p}B_{\rho +1}(y_{n, i})$ such that
\[
B_{\frac{R_{1}}{2}}(x_{0})\subset \mathbb{R}^{N}\backslash \cup_{j = 1}^{p}B_{\rho}(y_{n, i}).
\]
Next, for each $h, \eta >0$, let us consider
\[
\sigma_{h} = \frac{R_1}{2} + \frac{R_1}{2^{h+1}},\,\, \overline{\sigma}_{h} = \frac{\sigma_{h} + \sigma_{h+1}}{2}\quad\mbox{and}\quad K_{h} = \frac{\eta}{2}\Big(1-\frac{1}{2^{h+1}} \Big) \,\,\forall h = 0, 1, 2, ....
\]
Note that,
\[
\sigma_{h}\big\downarrow \frac{R_1}{2},\,\, K_{h}\big\uparrow \frac{\eta}{2}\quad\mbox{and}\quad \sigma_{h+1}<\overline{\sigma}_{h}< \sigma_{h}<1.
\]
In what follows, let us consider
\[
A_{n,K_{h}, \sigma_{h}}= \big\{x \in B_{\sigma_{h}}(x_0) \ : \ v_n(x) > K_{h} \big\}.
\]
For each $h = 0, 1, ...$, we fix
\[
J_{h, n} = \int_{A_{n,K_{h}, \sigma_{h}}}(( v_{n}-K_{h})_{+})^{\gamma^{*}}dx \quad \mbox{and}\quad \xi_{h}(x)= \xi\bigg(\frac{2^{h+1}}{R_1}\Big(\vert x-x_{0}\vert - \frac{R_1}{2}\Big) \bigg),
\]
where $\xi \in C^{1}(\mathbb{R})$ satisfies
\[
0\leq \xi\leq 1,\,\, \xi(t) = 1,\,\, \mbox{for}\,\,t\leq \frac{1}{2}\,\,\xi(t) = 0\,\,\mbox{for}\,\, t\geq \frac{3}{4}\,\,\mbox{and}\,\, |\xi^{'}|<c.
\]
Repeating the arguments explored in \cite[Lemma 3.5]{AS}, we can guarantee that
\[
J_{h+1, n} \leq CA^{h}J_{h, n}^{1 + \tau},
\]
where $C = C(N, \gamma, \gamma^{*}, R_1, \eta)$, $\tau = \frac{\gamma^{*}}{\gamma}-1$ and $A = 2^{\beta}$ for 
some $\beta$ sufficient large. We claim that there is $n_0 \in \mathbb{N}$ such that 
\begin{eqnarray}\label{J1}
	J_{0, n} \leq C^{\frac{1}{\tau}}A^{-\frac{1}{\tau^{2}}}, \forall n \geq n_0.
\end{eqnarray}
Indeed, note that
\begin{eqnarray*}
J_{0, n} = \int_{A_{n,K_{0}, \sigma_{0}}}(v_{n} - \frac{\eta}{2})_{+}^{\gamma^{*}}dx &\leq& \int_{A_{n,\frac{\eta}{2}, R_{1}}}(v_{n})_{+}^{\gamma^{*}}dx\\
&\leq& (\frac{\eta}{2})^{l^{*}}\int_{\mathbb{R}^{N}\backslash \cup_{j = 1}^{p}B_{\rho}(y_{n, i})}\Big(\frac{2v_{n}}{\eta}\Big)^{l^{*}}dx\\
& \leq& c_{1}\int_{\mathbb{R}^{N}\backslash \cup_{j = 1}^{p}B_{\rho}(y_{n, i})}\Phi_{*}(|v_{n}|)dx,
\end{eqnarray*}
where $c_1$ depends on $\eta$.  On the other hand, by Proposition \ref{TEmb} (see Appendix), there is $c_2>0$ independent of $\rho$ such that 
\[
\int_{\mathbb{R}^{N}\backslash \cup_{j = 1}^{p}B_{\rho}(y_{n, i})}\Phi_{*}(|v_{n}|)dx \leq c_{2}\Vert v_{n}\Vert_{W^{1, \Phi}(\mathbb{R}^{N}\backslash \cup_{j = 1}^{p}B_{\rho}(y_{n, i}))}.
\]
Hence, 
\[
J_{0, n}  \leq c_{3}\Vert v_{n}\Vert_{W^{1, \Phi}(\mathbb{R}^{N}\backslash \cup_{j = 1}^{p}B_{\rho}(y_{n, i}))},
\]
where $c_3>0$ depends on $\eta$. Now, using Claim \ref{fer1}, we can increase $\rho$, if necessary, of a way that 
$$
c_{3}\Vert v_{n}\Vert_{W^{1, \Phi}(\mathbb{R}^{N}\backslash \cup_{j = 1}^{p}B_{\rho}(y_{n, i}))}\leq C^{\frac{1}{\tau}}A^{-\frac{1}{\tau^{2}}},
$$
showing \eqref{J1}. Thus, by \cite[Lemma 4.7]{LU},
\[
\lim_{h \to +\infty}J_{h, n}= 0.
\]
On the other hand,
\[
\lim_{h \to +\infty}J_{h, n} =\lim_{h \to +\infty}\int_{A_{n,K_{h}, \sigma_{h}}}(( v_{n}-K_{h})_{+})^{\gamma^{*}}dx=\int_{A_{n,\frac{\eta}{2}, \frac{R_1}{2}}}(( v_{n}-\frac{\eta}{2})_{+})^{\gamma^{*}}dx,
\]
leading to
\[
v_{n}(z)\leq \frac{\eta}{2},\quad z \in  B_{\frac{R_1}{2}}(x_{0}),
\]
and so
\[
|v_{n}(z)|\leq \frac{\eta}{2},\quad z \in \mathbb{R}^{N}\backslash \cup_{j = 1}^{p}B_{\rho +1}(y_{n, i}),
\]
finishing the proof of the claim.

Hereafter, we consider the function $w_{n, i}(x)=v_{n}(x + y_{n, i})$. Note that it is a nonnegative and nontrivial solution of the problem
\[
- \Delta_{\Phi}w_{n, i} + V(\epsilon_{n}x + \epsilon_{n}y_{n, i})\phi(\vert w_{n, i}\vert)w_{n, i}  =  g(\epsilon_{n}x + \epsilon_{n}y_{n, i}, w_{n, i})\,\,\, \mbox{in} \,\, \mathbb{R}^{N}\eqno{(A_{\epsilon_{n}})}
\]
\begin{claim}\label{Ult1}
There exists $\delta>0$ such that $\Vert w_{n, i} \Vert_{\infty}\geq \delta$ for $n$ sufficient large.	
\end{claim}
In fact, if $\Vert w_{n, i} \Vert_{\infty}\to 0$,  \ref{f1} combined with $g$ gives 
\begin{eqnarray}\label{apend398}
	\frac{g(\epsilon_{n} x + \epsilon_{n}y_{n, i}, w_{n, i})}{\phi(|w_{n, i}|)w_{n, i}}\leq \frac{V_{0}}{2}\,\,\, \forall \,\, n\geq n_{0},
\end{eqnarray}
for some $n_{0} \in \mathbb{N}$. Now, \eqref{apend398} together with $J^{'}_{\epsilon_{n}}(w_{n, i})w_{n, i} = 0$ leads to
\begin{eqnarray*}
	\int_{\mathbb{R}^{N}}\phi(|\nabla w_{n, i}|)|\nabla w_{n, i}|^{2}dx + \int_{\mathbb{R}^{N}}V(\epsilon_{n} x + \epsilon_{n}y_{n, i})\phi(|w_{n, i}|)|w_{n, i}|^{2} dx =  0\,\,\, \forall \,\, n\geq n_{0},
\end{eqnarray*}
from where it follows that $\Vert w_{n, i} \Vert_{\epsilon_{n}} = 0$ for all $n\geq n_{0}$, which contradicts Lemma \ref{Nehari3}. 

\vspace{0.5 cm}

In the sequel, for $\eta < \delta$, the Claims \ref{afir989} and \ref{Ult1}  give
\[
\Vert w_{n, i} \Vert_{\infty, B_{(\rho+1)}(0)}\geq \delta,
\]
that is,
\[
\Vert v_{n} \Vert_{\infty, B_{\rho+1}(y_{n, i})}\geq \delta,\,\,\,\mbox{for all}\,\, i \in \Gamma.
\]
Finally, setting $u_{n}(x)= v_{n}\big(\frac{x}{\epsilon_{n}}\big)$ and $P_{n, i} = \epsilon_{n}y_{n, i}$, we get that $u_{n}$ is a solution of \eqref{PP3} verifying
\[
\Vert u_{n} \Vert_{\infty, B_{\epsilon_{n}(\rho+1)}(P_{n, i})}\geq \delta,\,\,\,\mbox{for all}\,\, i \in \Gamma.
\]
and
\[
\Vert u_{n} \Vert_{\infty, \mathbb{R}^{N}\backslash \cup_{i \in \Gamma}B_{\epsilon_{n}(\rho+1)}(P_{n, i})}\leq \Vert v_{n} \Vert_{\infty, \mathbb{R}^{N}\backslash B_{\rho+1}(y_{n, i})}\leq \eta \,\,\, \mbox{for all} \,\, n \geq n_{0},
\]
proving  the theorem.

\section{Appendix: New properties involving Orlicz-Sobolev spaces}
In this appendix, we will prove some results which were used in the present paper. Our first result is associated with an important property involving Orlicz-Sobolev spaces, which is well known for Sobolev spaces.  Here, we follows the same steps found in \cite[Theorem 3.2]{DT} (or \cite[Theorem 8.35]{Adams1}), however our proof can be applied for unbounded domains.
\begin{prop}\label{TEmb}
There exists $M^{*}>0$, which is independent of $\epsilon$, such that
\begin{eqnarray*}
\Vert u\Vert_{\Phi_{*}, \Omega_{\epsilon, i}}\leq M^{*} \Vert u \Vert_{\widetilde{X}_{\epsilon, i}}\,\,\, \mbox{for all}\,\, u \in \widetilde{X}_{\epsilon, i}.
\end{eqnarray*}
\end{prop}
\begin{pf}
In what follows, we define $\upsilon(t) = \big(\Phi_{*}(t)\big)^{1- \frac{1}{N}}$. Firstly, notice that
\begin{eqnarray}\label{apend30}
\Big|\frac{d}{dt}\upsilon(t)\Big|\leq \frac{N-1}{N}\widetilde{\Phi}^{-1}\big(\upsilon(t)^{\frac{N}{N-1}}\big)\,\,\,\mbox{for all}\,\, t>0.
\end{eqnarray}
For each $u \in \widetilde{X}_{\epsilon, i}\cap L^{\infty}(\Omega_{\epsilon, i})$ and $k>0$, the function $\nu:= \upsilon \circ \Big( \frac{|u|}{k}\Big) \in W^{1,1}(\Omega_{\epsilon, j})$ and
\[
\frac{\partial \nu (x)}{\partial x_{j}} =  \upsilon^{'}\Big(\frac{|u|}{k}(x)\Big)\frac{\mbox{sgn}u(x)}{k}\frac{\partial u (x)}{\partial x_{j}}.
\]
By \cite[Theorem 4.12]{Adams1},, once $\Omega_{_{\epsilon, j}}$ verifies the uniform cone condition for all $\epsilon>0$, we know that the constant associated with the embedding $ W^{1,1}(\Omega_{\epsilon, j}) \hookrightarrow L^{\frac{N}{N-1}}(\Omega_{\epsilon, j})$ does not depend on $\epsilon$, that is, there exists a positive constant $C$, which is independent of $u$ and $\epsilon$, such that
\[
\Vert \nu\Vert_{L^{\frac{N}{N-1}}(\Omega_{\epsilon, j})}\leq C\Big(\sum_{j=1}^{N}\Big\Vert \frac{\partial \nu}{\partial x_{j}}\Big\Vert_{ L^{1}(\Omega_{\epsilon, j})} + \Vert \nu\Vert_{ L^{1}(\Omega_{\epsilon, j})}\Big),
\]
or equivalently,
\begin{eqnarray*}
\Bigg[\int_{\Omega_{\epsilon, i}}\Phi_{*}\Big(\frac{|u|}{k}\Big)dx\Bigg]^{1- \frac{1}{N}}\leq \frac{C}{k}\sum_{j=1}^{N}\int_{\Omega_{\epsilon, i}} \Big|\upsilon^{'}\Big(\frac{|u|}{k}\Big)\frac{\partial u }{\partial x_{j}}\Big|dx + C\int_{\Omega_{\epsilon, i}} | \upsilon \Big( \frac{|u|}{k}\Big) | dx.
\end{eqnarray*}
Setting $k = \Vert u\Vert_{\Phi_{*}, \Omega_{\epsilon, i}}$, the Holder's inequality together with \eqref{apend30} yields
\begin{eqnarray}\label{apend33}
1\leq \frac{2C}{k}\frac{N-1}{N}\sum_{j=1}^{N} \Vert \widetilde{\Phi}^{-1}\big(\Phi_{*}\Big(\frac{|u|}{k}\Big)\big)\Vert_{\widetilde{\Phi}, \Omega_{\epsilon, i}}\Big\Vert \frac{\partial u }{\partial x_{j}}\Big\Vert_{\Phi, \Omega_{\epsilon, i}} + C\int_{\Omega_{\epsilon, i}} | \upsilon \Big( \frac{|u|}{k}\Big) | dx.\quad
\end{eqnarray}
Now, a direct computation leads to 
$$
\int_{\Omega_{\epsilon, i}} |\upsilon \Big( \frac{|u|}{k}\Big)| dx \leq  \frac{2}{k}\frac{N-1}{N}\big\Vert \widetilde{\Phi}^{-1}\big(\Phi_{*}(\frac{|u|}{k})\big)\big\Vert_{\widetilde{\Phi}, \Omega_{\epsilon, i}}\Vert u\Vert_{\Phi, \widetilde{\Omega}_{\epsilon, i}}.
$$
Since
$$
\Big\Vert \widetilde{\Phi}^{-1}\Big(\Phi_{*}\Big(\frac{|u|}{k}\Big)\Big)\Big\Vert_{\widetilde{\Phi}, \Omega_{\epsilon, i}}\leq 1,
$$
we get,
\begin{eqnarray}\label{apend35}
\int_{\Omega_{\epsilon, i}} |\upsilon \Big( \frac{|u|}{k}\Big)| dx &\leq& \frac{2}{k}\frac{N-1}{N}\Vert u\Vert_{\Phi, \widetilde{\Omega}_{\epsilon, i}}.
\end{eqnarray}
From \eqref{apend33}-\eqref{apend35}, 
\begin{eqnarray*}
1&\leq& \frac{2C}{k}\frac{N-1}{N} \Vert \nabla u \Vert_{\Phi, \widetilde{\Omega}_{\epsilon, i}}+ \frac{2}{k}\frac{N-1}{N}\Vert u\Vert_{\Phi, \widetilde{\Omega}_{\epsilon, i}}.\nonumber
\end{eqnarray*}
Hence, there exists $M_{*}>0$, independent of $\epsilon$ such that
\[
\Vert u\Vert_{\Phi_{*}, \Omega_{\epsilon, i}}\leq M^{*} \Vert u \Vert_{\widetilde{X}_{\epsilon, i}}\,\,\, \mbox{for all}\,\, u \in \widetilde{X}_{\epsilon, i}\cap L^{\infty}(\Omega_{\epsilon, i}),
\]
obtaining the desired result. 
\end{pf}

As a byproduct of the above proof, we have the following corollary 

\begin{corol} \label{NOVOCOR}
Let $(y_{n,i})$ the sequence obtained in (\ref{yni}). There is $C>0$, which is independent of $\rho$ and $n \in \mathbb{N}$, such that 
\begin{eqnarray*}
\int_{\mathbb{R}^{N}\backslash \cup_{j = 1}^{p}B_{\rho}(y_{n, i})}\Phi_{*}(|v|)dx \leq C \Vert v\Vert_{W^{1, \Phi}(\mathbb{R}^{N}\backslash \cup_{j = 1}^{p}B_{\rho}(y_{n, i}))},
\end{eqnarray*}
for all $v \in W^{1,\Phi}(\mathbb{R}^{N}\backslash \cup_{j = 1}^{p}B_{\rho}(y_{n, i}))$.
\end{corol}
\begin{pf} The corollary follows by repeating the same steps used in the proof Proposition  \ref{TEmb}. The main point that we would like to point out is the fact that the constant associated with the embedding 
$$ 
W^{1,1}(\mathbb{R}^{N}\backslash \cup_{j = 1}^{p}B_{\rho}(y_{n, i})) \hookrightarrow L^{\frac{N}{N-1}}(\mathbb{R}^{N}\backslash \cup_{j = 1}^{p}B_{\rho}(y_{n, i}))
$$
is also independent of $\rho$ and $n \in \mathbb{N}$, because $\Theta_{\rho,n,i}=\mathbb{R}^{N}\backslash \cup_{j = 1}^{p}B_{\rho}(y_{n, i})$ verifies the uniform cone condition for all $\rho>0$ and $n \in \mathbb{N}$.  
\end{pf}

The next result is also well known for Sobolev spaces, however for Orlicz-Sobolev spaces we do not know any reference. Here, we adapt some arguments found in \cite{AGJ}.
\begin{prop}\label{LEA3}
Let $\varrho>0$ and $\epsilon_{n}\in (0, +\infty)$ with $\epsilon_{n}\rightarrow 0$. Let $v_{n, i}\subset \widetilde{X}_{\epsilon_{n}}, i$ be a sequence and a constant $C_{0}>0$ such that
\[
\Vert v_{n, i}\Vert_{\widetilde{X}_{\epsilon_{n}, i}} \leq C_{0}\,\,\,\mbox{and}\,\,\,\lim_{n\rightarrow +\infty}\sup_{y \in \mathbb{R}^{N}}\int_{B_{\varrho}(y)\cap \Omega_{\epsilon_{n}, i}}\Phi(|v_{n, i}|)dx = 0.
\]
Then,
\[
\lim_{n\rightarrow +\infty}\int_{\Omega_{\epsilon_{n}, i}}B(|v_{n, i}|)dx = 0,
\]
for any N-function $B$ verifying $\Delta_{2}$-condition,
\[
\lim_{t\rightarrow 0}\frac{B(t)}{\Phi(t)} = 0\,\,\,\mbox{and}\,\,\, \lim_{|t|\rightarrow +\infty}\frac{B(t)}{\Phi_{*}(t)} = 0.
\]
\end{prop}
\begin{pf}
Firstly, note that given $\eta>0$ there exists $\kappa>0$ such that
\[
B(|v_{n, i}|)\leq \eta\Phi_{*}(|v_{n, i}|),\,\,\,\mbox{for}\,\, |v_{n, i}|\geq \kappa.
\]
As $(\|v_{n, i}\|_{\widetilde{X}_{\epsilon_{n}, i}})$ is bounded in $\mathbb{R}$, we have
\[
\int_{\Omega_{\epsilon_{n}, i}}B(|v_{n, i}|)dx \leq \eta C + \int_{\Omega_{\epsilon_{n}, i}\cap [|v_{n, i}|\leq \kappa]}B(|v_{n, i}|)dx
\]
which implies
\begin{eqnarray}\label{apend36}
\limsup_{n \to + \infty}\int_{\Omega_{\epsilon_{n}, i}}B(|v_{n, i}|)dx \leq \eta C + \limsup_{n \to + \infty}\int_{\Omega_{\epsilon_{n}, i}\cap [|v_{n, i}|\leq \kappa]}B(|v_{n, i}|)dx.
\end{eqnarray}
We will show that
\begin{eqnarray}\label{apend37}
\limsup_{n \to + \infty}\int_{\Omega_{\epsilon_{n}, i}\cap [|v_{n, i}|\leq \kappa]}B(|v_{n, i}|)dx = 0.
\end{eqnarray}
For this purpose, we consider for each $\zeta>0$ enough small, the function $\chi_{\zeta} \in C_{0}^{1}(\mathbb{R})$ given by 
\begin{align}
\chi_{\zeta}(s) =\left\{
\begin{array}
[c]{lcl}%
1,\,\,\,\mbox{if}\,\, |s|\le \kappa - \zeta,\\
a_{1}(s),\,\,\,\mbox{if}\,\,-(\kappa + \zeta) \leq s\leq -(\kappa - \zeta),\\
a_{2}(s),\,\,\,\mbox{if}\,\,\kappa - \zeta \leq s\leq \kappa + \zeta,\\
0,\,\,\,\mbox{if}\,\, |s|\geq \kappa + \zeta,
\end{array}
\right. \nonumber
\end{align}
where $a_{1}, a_{2}\in C^{1}\big(\mathbb{R};[0, 1]\big)$, $a_{1}$ is nondecreasing and $a_{2}$ is nonincreasing. Next, let us define the auxiliary function
\[
u_{n, i}(x) = \chi_{\zeta}(|v_{n, i}(x)|)v_{n, i}(x).
\]
Notice that
\begin{equation}\label{apend38}
\int_{\Omega_{\epsilon_{n}, i}}B(|u_{n, i}|)dx \geq \int_{\Omega_{\epsilon_{n}, i}\cap [|v_{n, i}|\leq \kappa - \zeta]}B(|v_{n, i}|)dx.
\end{equation}
Thereby,  (\ref{apend37}) follows by showing the limit below
\begin{equation} \label{apend39}
\limsup_{n \to + \infty}\int_{\Omega_{\epsilon_{n}, i}}B(|u_{n, i}|)dx = 0.
\end{equation}
In fact, gathering the above limit with \eqref{apend38}, we derive that 
\begin{eqnarray*}
\limsup_{n \to + \infty}\int_{\Omega_{\epsilon_{n}, i}\cap [|v_{n, i}|\leq \kappa - \zeta]}B(|v_{n, i}|)dx = 0.
\end{eqnarray*}
Since
\[
\int_{\Omega_{\epsilon_{n}, i}\cap [\kappa - \zeta \leq |v_{n, i}|\leq \kappa ]}B(|v_{n, i}|)dx = o_{n}(1)
\]
and
$$
\begin{array}{l}
\displaystyle \int_{\Omega_{\epsilon_{n}, i}\cap [|v_{n, i}|\leq \kappa]}B(|v_{n, i}|)dx = \int_{\Omega_{\epsilon_{n}, i}\cap [\kappa - \zeta \leq |v_{n, i}|\leq \kappa ]}B(|v_{n, i}|)dx \, + \\ 
\mbox{}\\
\hspace{4.5 cm} \displaystyle \int_{\Omega_{\epsilon_{n}, i}\cap [|v_{n, i}|\leq \kappa - \zeta]}B(|v_{n, i}|)dx,
\end{array}
$$
we deduce that
\[
\limsup_{n \to + \infty}\int_{\Omega_{\epsilon_{n}, i}\cap [|v_{n, i}|\leq \kappa]}B(|v_{n, i}|)dx=0,
\]
showing \eqref{apend37}. Now, by \eqref{apend36} and \eqref{apend37}, 
\begin{eqnarray*}
\limsup_{n \to + \infty}\int_{\Omega_{\epsilon_{n}, i}}B(|v_{n, i}|)dx \leq \eta C.
\end{eqnarray*}
By using that $\eta$ is arbitrary, it follows that 
\begin{eqnarray*}
\limsup_{n \to + \infty}\int_{\Omega_{\epsilon_{n}, i}}B(|v_{n, i}|)dx = 0,
\end{eqnarray*}
proving the proposition. Now, we observe that (\ref{apend39}) follows by repeating the same approach explored in \cite[Theorem 3.1]{AGJ}.

\end{pf}

\end{document}